\documentclass[12pt,reqno]{article}
\usepackage{amssymb}
\usepackage{amsmath}
\usepackage{amsthm}
\usepackage{color}
\usepackage{pb-diagram}
\usepackage{graphics,shortvrb}
\newcommand{\tw}[3]{{$#1$}${\,\scriptscriptstyle
{#2}}\atop\raise9pt\hbox{$\scriptstyle\tp$} ${$#3$}}
\newcommand{\st}[1]{\mbox{${\,\scriptscriptstyle
{#1}}\atop\raise5.5pt\hbox{$*$}$}}
\newcommand{\btr}{\raise1.2pt\hbox{$\scriptstyle\blacktriangleright$}\hspace{2pt}}
\newcommand{\id}{\mathrm{id}}

\newcommand{\tp}{\otimes}

\newcommand{\Lc}{\mathcal{L}}
\newcommand{\A}{\mathcal{A}}
\newcommand{\B}{\mathcal{B}}

\newcommand{\D}{\mathfrak{D}}
\newcommand{\M}{\mathrm{M}}
\newcommand{\R}{\mathrm{R}}

\newcommand{\Cc}{\mathcal{C}}
\newcommand{\Mc}{\mathcal{M}}

\newcommand{\C}{\mathbb{C}}

\newcommand{\Ha}{\mathcal{H}}
\newcommand{\Ru}{\mathcal{R}}

\renewcommand{\O}{\mathcal{O}}
\newcommand{\V}{V}
\newcommand{\U}{\mathcal{U}}
\newcommand{\F}{\mathcal{F}}
\newcommand{\ve}{\varepsilon}
\newcommand{\gm}{\gamma}

\newcommand{\ff}{\varphi}
\newcommand{\ot}{\otimes}

\newcommand{\la}{\lambda}

\newcommand{\Ab}{\mathrm{\bf A}}
\newcommand{\Bb}{\mathrm{\bf B}}

\newcommand{\End}{\mathrm{End}}

\newcommand{\Aut}{\mathrm{Aut}}

\newcommand{\Hom}{\mathrm{Hom}}
\newcommand{\Ob}{\mathrm{Ob}}

\newcommand{\tr}{\triangleright}
\newcommand{\tl}{\triangleleft}
\newcommand{\btl}{\mbox{\raise1.1pt\hbox{$\scriptstyle\blacktriangleleft$}}}

\newcommand{\g}{\mathfrak{g}}

\newcommand{\h}{\mathfrak{h}}

\newcommand{\n}{\mathfrak{n}}
\newcommand{\nn}{\nonumber}

\newcommand{\p}{\mathfrak{p}}
\renewcommand{\l}{\mathfrak{l}}
\renewcommand{\k}{\mathfrak{k}}
\renewcommand{\c}{\mathfrak{c}}
\newcommand{\K}{\mathcal{K}}

\newcommand{\MM}{\Lambda}

\newcommand{\si}{\sigma}
\newcommand{\al}{\alpha}
\newcommand{\bt}{\beta}
\newcommand{\be}{\begin{eqnarray}}
\newcommand{\ee}{\end{eqnarray}}

\newtheorem{thm}{Theorem}[section]
\newtheorem{propn}[thm]{Proposition}
\newtheorem{lemma}[thm]{Lemma}
\newtheorem{corollary}[thm]{Corollary}
\newtheorem{conjecture}[thm]{Conjecture}

\theoremstyle{definition}
\newtheorem{remark}[thm]{Remark}
\newtheorem{definition}[thm]{Definition}
\newtheorem{example}[thm]{Example}

\newcommand{\select}[1]{\textcolor{red}{\bf\em #1}}

\begin{document}
\title{Dynamical Yang-Baxter equation and quantum vector bundles\footnote{
The research is supported in part
by the Israel Academy of Sciences grant no. 8007/99-03, the Emmy Noether
Research Institute for Mathematics,
 the Minerva Foundation of Germany, the Excellency Center "Group
Theoretic Methods in the study of Algebraic Varieties"  of the Israel
Science foundation, and by Russian Foundation for Basic Research  grant no.
03-01-00593.}}
\author{J. Donin
\hspace{3pt} and A. Mudrov}
\date{}
\maketitle
\begin{center}
{Department of Mathematics, Bar Ilan University, 52900 Ramat Gan,
Israel.}
\end{center}
\begin{abstract}
We develop a categorical approach to the dynamical Yang-Baxter equation
(DYBE)
for arbitrary Hopf algebras. In particular, we introduce the notion of a
dynamical
extension of a monoidal category, which provides a natural environment
for quantum dynamical R-matrices, dynamical twists, {\em etc}. In this
context, we define
dynamical associative algebras and show that such algebras give
quantizations of vector
bundles on coadjoint orbits.
We build a dynamical twist for any pair of a reductive Lie algebra and
their
Levi subalgebra. Using this twist, we obtain an equivariant star product
quantization of vector bundles on semisimple coadjoint orbits of reductive
Lie groups.
\end{abstract}
{\small \underline{Key words}: Dynamical Yang-Baxter equation, dynamical
categories, quantum vector bundles.\\
\underline{AMS classification codes}: 17B37, 81R50.}
\tableofcontents
\section{Introduction}

The quantum dynamical Yang-Baxter equation (DYBE)  appeared in the
mathematical physics literature, \cite{GN,AF,Fad,F,ABB},
in connection with integrable models of conformal field theories.
The classical DYBE was first considered
in \cite{BDFh}, rediscovered in \cite{F}, and systematically studied  by
Etingof, Schiffmann, and Varchenko, \cite{EV1,ES2,S}.
For a guide in the DYBE theory and an extended bibliography the reader is
referred to the lecture course
\cite{ES1}.

The theory of DYBE over the Cartan subalgebra in a simple Lie algebra has
been developed in detail.
Classical dynamical r-matrices were classified in \cite{EV1} and their
explicit quantization
built in \cite{EV2,EV3,ESS}.
Concerning the  classical DYBE over an arbitrary (non-commutative) base,
much is known about classification
of its solutions and there are numerous explicit examples,
\cite{AM,ES2,Fh,S,Xu2}.
At the same time, there is no generally accepted
definition of quantum DYBE over a non-commutative Lie algebra or, say, over
an arbitrary Hopf
algebra. A generalization of the quantum DYBE for several particular cases
was proposed in \cite{Xu2} and  \cite{EE1}.
Such a generalization was motivated by a relation between DYBE and
the star product, \cite{Xu1,Xu2}. An open question is an
interpretation of the quantum DYBE of  \cite{Xu2} and  \cite{EE1}
from a categorical point of view.

Another interesting question is
a relation of DYBE to the equivariant quantization.
It is observed by Lu, \cite{Lu1}, that the list of classical r-matrices
over the Cartan subalgebra of a simple
Lie algebra is in intriguing correspondence with the list of
Poisson-Lie structures on its maximal coadjoint orbits.
However, the precise relation between quantum dynamical R-matrices
and the equivariant quantization  has not been clarified.

The purpose of the present paper is to develop the theory of DYBE over an
arbitrary
Hopf algebra and relate it to equivariant quantization of vector bundles.
Firstly, we generalize the classical dynamical Yang-Baxter equation for
any Lie bialgebra $\h$ extending the concept of base manifold, which is the
dual space $\h^*$
in the standard approach.
Secondly, we build dynamical extensions of monoidal categories and
define the quantum dynamical R-matrix over arbitrary base.
Our third result is a construction of dynamical twist for Levi subalgebras
in
a reductive Lie algebra. Finally, we introduce a notion of dynamical
associative algebras
as algebras in dynamical categories. We relate them to equivariant
quantization
of vector bundles. As an application, we construct an equivariant star
product
quantization of vector bundles  (including function algebras) on
semisimple
coadjoint orbits of reductive Lie groups.

It turns out that there is a general procedure of "dynamical extension",
$\bar \O$,
of every monoidal category $\O$ over a base $\B$, which is an $\O$-module
category.
This new category has the same objects as $\O$ but more morphisms.
The objects are considered as functors from $\B$ to $\B$ by the tensor
product action.
Morphisms are natural transformations between these functors.
This category admits a tensor product making it a monoidal category with
$\O$ being a subcategory. One can consider the standard notions
as algebras, twists, and R-matrices relative to $\bar \O$.
In terms of the original category $\O$, they  satisfy "shifted"
axioms, like shifted associativity, shifted cocycle condition,
shifted or dynamical Yang-Baxter equation.

The construction of dynamical extension admits various formulations.
One of them uses the so called base algebras, which
are commutative algebras in the Yetter-Drinfeld categories.
From the algebraic point of view, a Yetter-Drinfeld category is
a category of modules over the double $D(\Ha)$ of a Hopf algebra $\Ha$.
In the quasi-classical limit, the base algebras are function algebras
on the so-called Poisson base manifolds.
A Poisson base manifold $L$ is endowed with an action of the double
$D(\h)$
of the Lie bialgebra $\h$, the classical analog of $\Ha$.
The Poisson structure on $L$ is induced by the canonical $r$-matrix of
the double.

The category of $\Ha$-modules can be dynamically extended over  the dual
Hopf algebra $\Ha^*$.
This approach is convenient for definition of dynamical associative
algebras.
A dynamical associative algebra is equipped with an equivariant family of
binary operations (multiplications)
depending on elements of $\Ha^*$.
This family satisfies a "shifted" associativity condition.
We show that the dynamical associative algebras give vector bundles
on quantum spaces.

In this paper we consider vector bundles on coadjoint orbits.
In the classical situation, the function algebra on a homogeneous space
is a subalgebra in the function algebra on the group.
In general, the  quantized function algebra on a homogeneous space cannot
be realized as
a subalgebra in a quantized function algebra on the group.
For example, in the case of semisimple coadjoint orbits,
such a realization exists only for symmetric or bisymmetric orbits,
\cite{DGS1,DM1}.
Nevertheless,  the quantization of the function algebra on the group
to a \select{dynamical associative algebra} contains quantum orbits as
(associative)
subalgebras. Moreover, a dynamical quantization on the group
quantizes the algebra of sections of homogeneous vector bundles on orbits.
Such quantizations  are parameterized by group-like elements of $\Ha^*$.

A way of constructing (quantum) dynamical R-matrices and dynamical
associative
algebras is by twists in dynamical categories. We build such twists for
Levi subalgebras in simple Lie algebras, using
generalized Verma modules. This gives a construction of star product
on the semisimple orbits.

The paper is organized as follows. In Section \ref{sDRM} we recall
basic definitions concerning DYBE and  compatible
star product of \cite{Xu2}.

Section  \ref{sDYA} presents generalizations of
DYBE using the concepts of base algebras and base manifolds.

Section \ref{sDC} is devoted to various formulations of dynamical
categories,
therein we study dynamical associative algebras.

In Section \ref{sDCRA} we study objects that are interesting for
applications:
dynamical twists and dynamical R-matrices.
We consider various types of dynamical categories and give
expressions of dynamical twists and R-matrices in terms
of the original category.

In Section \ref{sDCC} we suggest  a method of constructing dynamical
twists.
The method is based on  a notion of dynamical adjoint functors.
We build such functors using generalized Verma modules
corresponding to Levi subalgebras in (quantum) universal enveloping
algebra of simple Lie algebras.

In Section \ref{sQVB} study relations between  quantization of vector
bundles
and dynamical  associative algebras in a purely algebraic setting.

In Section \ref{sVBCO} we give a detailed consideration to the dynamical
associative algebra which is a quantized function algebra on
the Lie group $G$. We relate this algebra to quantum vector bundles on
coadjoint semisimple orbits of $G$.

Note that the equivariant star product on function algebras on coadjoint
orbits
was also constructed in the papers \cite{AL} and \cite{KMST} appeared
after
the first version of this paper.
Our method of building dynamical twists is developed for a more general
case in \cite{EE2}.

\vspace{0.3cm}
\noindent
{\bf \large Acknowledgements}.
\hspace{10pt}
We are grateful to J. Bernstein, V. Ostapenko, and S. Shnider for
stimulating discussions within the "Quantum groups" seminar
at the Department of Mathematics, Bar Ilan University.
We appreciate useful remarks by M. Gorelik,  V. Hinich, and A. Joseph
during a talk at the Weizmann Institute. Our special thanks to P. Etingof
for his comments on various aspects of the subject.
\section{Dynamical r-matrix and compatible star product}
\label{sDRM}
\subsection{Classical dynamical Yang-Baxter equation}
In this section we recall basic definitions concerning the
dynamical Yang-Baxter equation. Let $\g$ be a Lie algebra and $\h$
its Lie subalgebra. The dual space $\h^*$ is considered as an
$\h$-module with respect to the coadjoint action. Let
$\{h_i\}\subset \h$ be a basis and $\{\la^i\}\subset \h^*$ its
dual.
\begin{definition}[\cite{F,EV1}]
\label{def_0} A classical dynamical r-matrix over the base $\h$ is
an $\h$-equiva\-riant meromorphic function $r\colon\h^*\to \g\tp
\g$ satisfying
\begin{enumerate}
\item the normal condition: the sum $r(\la)+r_{21}(\la)$ is
$\g$-invariant, \item  the classical dynamical Yang-Baxter
equation (DYBE):
\be
\label{cDYBE0}
\sum_{i} \frac{\partial r_{23}}{\partial \la^i}h^{(1)}_i-
\frac{\partial r_{13}}{\partial \la^i}h^{(2)}_i+\frac{\partial
r_{12}}{\partial \la^i}h^{(3)}_i
 = [r_{12},r_{13}]+[r_{13},r_{23}]+[r_{12},r_{23}].
\ee
\end{enumerate}
\end{definition}
\noindent A constant dynamical $r$-matrix is a solution to the
ordinary Yang-Baxter equation. It follows that the sum
$r(\la)+r_{21}(\la)$ does not depend on $\la$, \cite{ES2}. If it
is identically zero, the r-matrix is called triangular.
\subsection{Quantum dynamical Yang-Baxter equation over a commutative
base}
Suppose $\h$ is a commutative Lie algebra and $V$ is a semisimple
$\h$-module. Given a family $\Omega(\la)$, $\la\in\h^*$, of linear
operators on $V^{\tp 3}$, let us denote by $\Omega(\la+t h^{(1)})$
the family of operators on  $V^{\tp 3}$ acting by $v_1\tp v_2\tp
v_3\mapsto\Omega(\la+t\; \mathrm{wt}(v_1))(v_1\tp v_2\tp v_3)$
where $\mathrm{wt}(v)$ stands for the weight of $v\in V$ with
respect to $\h$ and $t$ is a formal parameter. The
operators $\Omega(\la+t h^{(i)})$, $i=2,3$, are defined similarly.

\begin{definition} Let $\h$ be a commutative Lie subalgebra of a
Lie algebra $\g$. Let $R(\la)$ be an $\h$-equivariant meromorphic
function $\h^*\to \U(\g)^{\tp 2}$ (we consider $\h^*$ equipped
with the coadjoint and $\U(\g)$ with adjoint action of $\h$). Then
$\Ru(\la)$ is called  (universal) quantum dynamical R-matrix if
it satisfies the quantum dynamical Yang-Baxter equation (QDYBE)
\be \label{DYBE}
\Ru_{12}(\la)\Ru_{13}(\la+t h^{(2)})\Ru_{23}(\la)
=
\Ru_{23}(\la+ t h^{(1)})\Ru_{13}(\la)\Ru_{12}(\la+t h^{(3)}).
\ee
\end{definition}

Assuming  $\Ru(\la)=1\tp 1 + t \:r(\la)+O(t^2)$, the element
$r(\la)$ satisfies equation (\ref{cDYBE0}), i.e. equation
(\ref{cDYBE0}) is the quasi-classical limit of equation
(\ref{DYBE}). In this case $R(\la)$ is called quantization of
$r(\la)$. The  problem of quantizing classical DYBE has been solved for
$\g$ a complex semisimple Lie algebra and $\h$ its reductive
commutative subalgebra, \cite{ESS}. As to the case of general
$\h$, there is no generally accepted concept of what should be
taken as the quantum DYBE. In the next subsection we render a
construction of \cite{Xu2} suggesting a version of quantum DYBE as
a quantization ansatz for triangular dynamical r-matrices.
This will be
the starting point for our study.
\subsection{Compatible star product}
Let $\g$ be a complex  Lie algebra and $G$ the corresponding
connected Lie group. Let $\h$ be a Lie subalgebra  in $\g$. Denote
by $\vec\xi$ the left invariant vector field on $G$ induced by
$\xi\in \g$ via the right regular action. Let $\pi_{\h^*}$ denote
the Poisson-Lie bracket on $\h^*$.
\begin{thm}[\cite{Xu2}]
\label{Xu} A smooth function $r\colon\h^*\to \wedge^2 \g$ is a
triangular dynamical r-matrix if and only if the bivector field
\be
\label{DPB}
\pi_{\h^*}+\sum_i\frac{\partial }{\partial \la_i}\wedge \vec{h_i} +
\vec{r}(\la)
\ee
is a Poisson structure on $\h^*\times G$.
\end{thm}
\noindent
Thus, the bivector field $\vec r(\la)$ on $G$ is a
"part" of a special Poisson bracket on a bigger space, $\h^*\times
G$. Xu proposed to look at a star product on $\h^*\times G$  of
special form, as a quantization of (\ref{DPB}). Let
$\h_t:=\h[[t]]$ be a Lie algebra over $\C[[t]]$ with the Lie
bracket $[x,y]_t:=t[x,y]$ for $x,y\in\h$. The algebra $\U(\h_t)$
can be considered as a deformation quantization of the polynomial
algebra on $\h^*$. It is known that this quantization can be
presented as a star product on $\h^*$ by the PBW map
$S(\h)[[t]]\to\U(\h_t)$, where elements of the symmetric algebra
$S(\h)$ are identified with polynomial functions on $\h^*$. We
call this star product the \select{PBW star product}.

\begin{definition}[\cite{Xu2}]
\label{defCSP}
A star product $*_t$ on $\h^*\times G$  is called
\select{compatible} if
\begin{enumerate}
\item
when restricted to $C^\infty(\h^*)$, it coincides with the
PBW star product; \item for $f\in C^\infty(G)$ and $g\in
C^\infty(\h^*)$
\be
(f*_t g)(\la,x):=f(x)g(\la),\quad (g*_t f)(\la,x)
:=
\sum_{k=0}^\infty \frac{t^k}{k!} \frac{\partial^k
g(\la)}{\partial\la^{i_1}
\ldots\partial\la^{i_k}} \vec h_{i_1}\ldots \vec h_{i_k} f(x);
\label{eq1}
\ee
\item
for $f,g\in C^\infty(G)$
\be
\label{eqforF}
(f*_t g)(\la,x):=\vec \F(\la) (f,g)(x),
\ee
where $\F(\la)$ is a smooth function $\F\colon \h^*\to
\U(\g)\tp\U(\g)[[h]]$
such that $\F=1\tp 1+ \frac{t}{2}\:r(\la)+ O(t^2)$.
\end{enumerate}
\end{definition}
\noindent For this star product to be associative, $\F$ should
satisfy a certain condition called shifted cocycle condition.

Also, Xu proposed a generalization of the quantum DYBE (\ref{DYBE}) for
arbitrary Lie algebra $\h$ in the form
\be
\label{XuYBE}
\Ru_{12}(\la)*_t \Ru_{13}\bigl(\la+t h^{(2)}\bigr)*_t \Ru_{23}(\la)
=
\Ru_{23}\bigl(\la+t h^{(1)}\bigr)*_t \Ru_{13}*_t \Ru_{12}\bigl(\la+t
h^{(3)}\bigr),
\ee
where $\Ru$ is an equivariant function $\h^*\to \U(\g)\tp \U(\g)$,
 and the subscripts mark
the tensor components in $\U^{\tp 3}(\g)$. Notation $f(\la+t h)$
for $f\in C^\infty (\h^*)$ means
\be
f(\la+t h):=\sum_{k=0}^\infty \frac{t^k}{k!}\frac{\partial^k
f(\la)}{\partial\la^{i_1}\ldots\partial\la^{i_k}}\; h_{i_1}\ldots
h_{i_k}.
\label{eqPBWco}
\ee
Here $\{h_i\}\subset \h$ and
$\{\la^i\}\subset \h^*$ are dual bases; the superscript of $h^{(i)}$,
$i=1,2,3$, in (\ref{XuYBE}) means that $\h$ is embedded  in the
$i$-th component of $\U^{\tp 3}(\g)$. A version of compatible star
product with non-trivial associativity was proposed in \cite{EE1}
for the quantization of the Alekseev-Meinrenken dynamical r-matrix
over $\h=\g$, \cite{AM}.

The compatible star product of \cite{Xu2} is defined on smooth
functions on $\h^*\times G$. When restricted to polynomial
functions on $\h^*$, it gives the multiplication in the universal
enveloping algebra $\U(\h)$. Formula (\ref{eq1}) expresses the
product of elements from $\U(\h)$ and $C^\infty(G)$ through the
comultiplication in  $\U(\h)$ and the action of $\U(\h)$ on
$C^\infty(G)$. It seems natural to replace  $\U(\h)$ with an
arbitrary Hopf algebra $\Ha$ and $C^\infty(G)$ with a right
$\Ha$-module $\A$. However, the bidifferential operator
$\F(\la)$ in (\ref{eqforF}) may be a meromorphic or even a formal
function in $\la\in\h^*$. This requires to consider
 appropriate extensions of $\U(\h)$, which may no
longer  be Hopf algebras. On the other hand, there is a class of
\select{admissible} algebras which are close, in a sense, to the Hopf
ones. Those are commutative algebras in the so-called
Yetter-Drinfeld category of $\Ha$-modules and $\Ha$-comodules,
which are, roughly speaking, modules over the double of $\Ha$. We
will define a \select{dynamical extension} of the monoidal category
of $\Ha$-modules over an admissible algebra, where the notions of
compatible star products, dynamical Yang-Baxter equations {\em etc},
acquire a natural interpretation. Depending on a particular
choice of an admissible algebra, we come to different
quasi-classical limits of quantum dynamical objects. Also, it
appears useful (and often technically simpler) to consider a
"dual" version of the dynamical extension, for example, a
dynamical extension of the monoidal category of $\Ha^*$-comodules.
In this way we obtain a "linearization" of the theory; in
particular, smooth or meromorphic functions on $\h^*$ become
linear functions on $\U(\h)^*$. Moreover, it will be useful to
introduce the notion of dynamical extension of an arbitrary
monoidal category, defined without involving any Hopf algebra.
Below we present all the formulations.
\section{Generalizations of dynamical Yang-Baxter equations}
\label{sDYA}
\subsection{Base algebras}
\label{ssBA}
In this subsection we define two objects of our primary concern: a
base algebra  $\Lc$ and a dynamical associative algebra over
$\Lc$.

By $k$ we mean a commutative ring over a field of zero
characteristic. The reader may think of it as $\C$ or $\C[[t]]$,
the ring of formal series in $t$. Given a Hopf algebra $\Ha$
over $k$ we denote  the multiplication, comultiplication,  counit,
and antipode by $m$, $\Delta$, $\ve$, and $\gamma$. We use the
standard Sweedler notation for the comultiplication in Hopf
algebras: $\Delta(x)=x^{(1)}\tp x^{(2)}$. In the same fashion we
denote the $\Ha$-coaction on a right comodule $A$:
$\delta(a)=a^{[0]}\tp a^{(1)}$, where the square brackets label
the $A$-component and the parentheses mark that belonging to
$\Ha$. The Hopf algebra with the opposite multiplication will be
denoted by $\Ha_{op}$ while with the opposite comultiplication  by
$\Ha^{op}$.

The Hopf algebra $\Ha$ is considered as a left module
over itself with respect to the adjoint action
\be
\label{eq_ad}
x\tp a \mapsto x^{(1)} a \gm(x^{(2)});
\ee
then the multiplication
in $\Ha$ is equivariant. It is a standard fact that for any  left
$\Ha$-module $A$ the map $\Ha \tp A\to  A\tp \Ha$,
$h\tp a\mapsto h^{(1)}\tr a \tp h^{(2)}$, is $\Ha$-equivariant.

Recall that an algebra and $\Ha$-module $\A$ is a module algebra
if the multiplication in $\A$ is $\Ha$-equivariant. Dually, an
algebra and $\Ha$-comodule $\A$ is a comodule algebra if the
coaction $\A\to\Ha\ot\A$ is a homomorphism of algebras.

\begin{definition}[Base algebras]
\label{def_BA}
A left  $\Ha$-module and left $\Ha$-comodule
algebra $\Lc$ is called \select{base algebra} over $\Ha$ if the coaction
$\delta\colon\Lc\to \Ha\tp \Lc$ satisfies the condition
\be
\label{eqDYM}
\bigl\{ x^{(1)}\tr \ell\bigr\}^{(1)} x^{(2)}\tp \bigl\{ x^{(1)}\tr
\ell\bigr\}^{[2]}
=
 x^{(1)}\ell^{(1)}   \tp x^{(2)}\tr \ell^{[2]}
\ee
for all $x\in \Ha$ and $\ell\in \Lc$, and the condition
\be
\label{eqQCDYA}
\ell_1 \ell_2= \bigl(\ell_1^{(1)} \tr \ell_2\bigr) \:\ell_1^{[2]},
\ee
for all $\ell_1,\ell_2\in \Lc$.
\end{definition}
The coaction $\delta$ defines a permutation $\tau_A \colon \Lc\tp A\to A\tp
\Lc$
with every $\Ha$-module $A$:
\be
\label{permutation}
\tau_A(\ell\tp a):= \ell^{(1)}\tr a\tp \ell^{[2]},\quad \ell\tp a\in \Lc\tp
A.
\ee
Condition (\ref{eqDYM}) ensures that this permutation is
$\Ha$-equivariant.
Condition (\ref{eqQCDYA}) means that the multiplication in $\Lc$ is
$\tau_\Lc$-commutative.
\begin{remark}
A base algebra is a
commutative algebra in the braided category of Yetter-Drinfeld modules.
From the purely algebraic point of view, Yetter-Drinfeld modules
are modules over the double Hopf algebra
$D(\Ha)$. The left $\Ha$-coaction induces a left $\Ha^*_{op}$-action.
Together with the $\Ha$-action, the $\Ha^*_{op}$-action gives a
$D(\Ha)$-action.
In our theory, an $\Ha$-base algebra plays the same role as
the $\U(\h)$-module algebra of functions on $\h^*$ in the theory of DYBE
over
a commutative base.

One can also introduce the dual notion of a base coalgebra
as a comodule over $D(\Ha)$. We
will use
 $\Ha^*$, a dual to the Hopf algebra $\Ha$, as an example of
 such a base coalgebra.
\end{remark}

\begin{example}
The algebra $\Ha$ itself is a base algebra over $\Ha$ with respect to
the left adjoint action and the coproduct $\Delta$
considered as the left regular $\Ha$-coaction. Conditions
(\ref{eqDYM}) and (\ref{eqQCDYA}) are checked directly.
\end{example}
\begin{example}
\label{basereduction}
Suppose that $\Ha$ is the tensor product of two Hopf algebras,
$\Ha=\Ha_0\tp \Ha_1$.
Then both $\Ha_0$ and $\Ha_1$ are  naturally base algebras over $\Ha$.
The $\Ha$ action on $\Ha_i$ is the adjoint action restricted
to $\Ha_i$.
The $\Ha$-coaction on $\Ha_i$ is the coproduct coaction considered
as a map with values in $\Ha_i\tp \Ha_i\subset \Ha\tp \Ha_i$.
\end{example}
\begin{example}[PBW star product]
\label{PBW*pr}
Consider the algebra $C^\infty(\h^*)[[t]]$ from Definition \ref{defCSP}
equipped with the PBW star product. It is obviously a left
$\U(\h)$-module algebra, and formula (\ref{eqPBWco}) defines a
coproduct $C^\infty(\h^*)[[t]]\to \U(\h)\tp C^\infty(\h^*)[[t]]$ (the
completed
tensor product). It is
straightforward to check that  $C^\infty(\h^*)[[t]]$ is a  base algebra
over $\U(\h)$. The algebra  $C^\infty(\h^*)[[t]]$ is an extension of
$\U(\h_t)$, which is realized as the subalgebra in $\U(\h)[[t]]$
generated by $t\h$. The algebra $\U(\h_t)$ is a Hopf one, hence it
is a base algebra over itself. At the same time, it is a base
algebra over $\U(\h)[[t]]$. Indeed, it is invariant under the adjoint
$\U(\h)$-action, and it is a left $\U(\h)$-comodule under the map
$(\varphi_t\tp \id)\circ \Delta$, where $\Delta$ is the coproduct
in $\U(\h_t)$ and $\varphi_t$ the natural embedding of $\U(\h_t)$
in $\U(\h)[[t]]$.
\end{example}

\begin{propn}
\label{exLcom}
Suppose that $\Ha$ is quasitriangular Hopf algebra, with the
universal R-matrix $\Ru$. Let $\Lc$ be a quasi-commutative $\Ha$-module
algebra, i.e. obeying $(\Ru_2\tr \ell_2)(\Ru_1\tr \ell_1) = \ell_1\tp
\ell_2$
for all $\ell_1,\ell_2\in \Lc$.
Then $\Lc$ is an $\Ha$-base algebra, with the left $\Ha$-coaction
\be
\label{R-coaction}
\delta(\ell):=\Ru_2\tp \Ru_1\tr \ell, \quad \ell\in \Lc.
\ee
\end{propn}
\begin{proof}
The condition (\ref{eqQCDYA}) is satisfied by construction.
The equality $(\Delta\tp \id)(\Ru)=\Ru_{13}\Ru_{23}$ implies that
the map (\ref{R-coaction}) is an algebra homomorphism.
The map (\ref{R-coaction}) makes $\Lc$ a left $\Ha$-comodule,
because of $(\id\tp \Delta)(\Ru)=\Ru_{13}\Ru_{12}$.
The condition (\ref{eqDYM}) holds by virtue of
$\Ru\Delta(h)=\Delta^{op}(h)\Ru$ for every $h\in \Ha$.
\end{proof}
\begin{corollary}
\label{corR-coaction}
Within the hypothesis of  Proposition \ref{exLcom}, suppose
that $\Ru\in \Ha\tp \K\subset \Ha\tp \Ha$, where $\K$ is  a Hopf subalgebra
in $\Ha$
Then $\Lc$ is endowed with a structure of $\K$-base algebra,
with the $\K$-coaction (\ref{R-coaction}).
\end{corollary}
\begin{proof}
The $\Ha$-coaction (\ref{R-coaction}) is, in fact, an $\K$-coaction.
Now the statement immediately follows from Proposition \ref{exLcom}.
\end{proof}

Remark that an $\Ru$-commutative algebra $\Lc$ is commutative with respect
to the element $\Ru^{-1}_{21}$,
which is also a universal R-matrix for $\Ha$. Thus $\Lc$ has
two $\Ha$-base algebra structures, and they are different in general.
In particular, an $\Ha$-base algebra has two different $\D(\Ha)$-base
algebra structures.

\begin{example}[The FRT algebras]
\label{FRT}
The FRT-dual Hopf algebra $\Ha^*$, \cite{FRT}, of a quasitriangular Hopf
algebra $\Ha$ is a quasi-commutative $\Ha\tp \Ha_{op}$-algebra.
Therefore it has two structures of $\Ha\tp \Ha_{op}$-base algebras.
\end{example}
\begin{example}[Reflection equation algebras]
\label{exREa}
Recall that a twist of a Hopf algebra $\Ha$ is a
Hopf algebra with the same multiplication and the new
comultiplication $\tilde\Delta(x):=\F^{-1}\Delta(x) \F$, where the
element $\F$ called a twisting cocycle satisfies certain conditions, see
\cite{Dr3}. For every quasitriangular Hopf algebra $\Ha$ with the
R-matrix $\Ru$, there is the twist \tw{\Ha}{\Ru}{\Ha} of its tensor
square, \cite{RS}. It is obtained by applying the twisting cocycle
$\Ru_{23}\in
(\Ha\tp \Ha)\tp (\Ha\tp \Ha)$ to the
comultiplication in $\Ha\ot\Ha$. The twisted tensor square is a
quasitriangular Hopf
algebra with the R-matrix
\be
\label{TTPR}
\Ru':=\Ru_{14}^{-}\Ru_{13}^{-}\Ru_{24}^+\Ru_{23}^+
\in (\mbox{\tw{\Ha}{\Ru}{\Ha}})\tp (\mbox{\tw{\Ha}{\Ru}{\Ha}}),
\ee
where $\Ru^+ := \Ru$ and $\Ru^- := \Ru^{-1}_{21}$.
Recall that $\Ha$ is a Hopf subalgebra in \tw{\Ha}{\Ru}{\Ha}
through  the embedding $\Delta\colon \Ha\to \Ha\tp \Ha$.
Observe that the R-matrix (\ref{TTPR})  can be presented
as $\Ru'=(\Ru^-_1\tp \Ru^+_1)\tp \Delta(\Ru^-_2 )\Delta(\Ru^+_2 )$.
In other words, its right tensor component belongs to
$\Delta(\Ha)\subset \mbox{\tw{\Ha}{\Ru}{\Ha}}$.
Applying the argument form Corollary \ref{corR-coaction}
to $\K=\Delta(H)$, we come to the following proposition.
\begin{propn}
\label{qcHHA}
A quasi-commutative \tw{\Ha}{\Ru}{\Ha}-module algebra is a base
algebra over $\Ha$.
\end{propn}

The reflection equation algebra associated with a finite dimensional
representation of $\Ha$, \cite{KSkl,KS}, is a quasi-commutative
\tw{\Ha}{\Ru}{\Ha}-algebra,
\cite{DM3}. As a corollary of Proposition \ref{qcHHA}, we obtain
that the reflection equation algebra is an $\Ha$-base algebra.

More examples of base algebras are obtained by quantizing
Poisson base algebras, see examples of Section \ref{ssPBAPBM}, according to
Theorem \ref{BAquantization}.
\end{example}
\subsection{Dynamical associative algebras}
\label{ssDAA}
Let $\Lc$ be a base algebra over a Hopf algebra $\Ha$.
\begin{definition}
\label{defDA} A left $\Ha$-module $\A$ is called \select{dynamical
associative  algebra} over the base algebra $\Lc$ if it is equipped
with an $\Ha$-equivariant bilinear map $\divideontimes\colon\A\tp
\A\to \A\tp \Lc$ such that the following diagram is commutative:
\be
&
\begin{diagram}
   \dgARROWLENGTH=0.5\dgARROWLENGTH
\node{\A \tp \Lc\tp \A}\arrow{e,t}{\id \tp\tau_\A}\node{\A \tp \A \tp \Lc}
\arrow{e,t}{\divideontimes\tp \id}\node{\A \tp \Lc \tp
\Lc}\arrow{e,t}{\id\tp\mathrm{m}}\node{\A \tp \Lc}
\arrow{s,r,!}{\parallel}
\\
\node{\A\tp\A\tp \A }\arrow{e,t}{\id\tp \divideontimes
}\arrow{n,l}{\divideontimes\tp \id}
\node{\A \tp \A \tp \Lc}\arrow{e,t}{\divideontimes\tp \id}\node{\A \tp \Lc
\tp \Lc}\arrow{e,t}{\id\tp\mathrm{m}}\node{\A \tp \Lc}
\end{diagram}
\label{shifted_assoc0}
\ee
Here $\mathrm{m}$ stands for the multiplication in $\Lc$
and the permutation $\tau_\A$ is defined by (\ref{permutation}).
\end{definition}

An example of dynamical associative algebra is the function algebra on a
group
$G$ twisted by the dynamical twist from \cite{Xu2}. It defines
the compatible star product in the sense of Definition \ref{defCSP};
it turns out that the multiplication $\divideontimes$ in a dynamical
associative algebra over an
arbitrary base can be extended
to an ordinary associative multiplication in a bigger algebra, according to
the following proposition.
\begin{propn}
\label{propCSP}
 Let $\A$ be a left $\Ha$-module equipped
with an equivariant map $\divideontimes\colon\A\tp \A\to \A\tp
\Lc$. Then $\A$ is a dynamical associative  algebra with respect
to $\divideontimes$ if and only if the operation
\be
(\A\tp \Lc)\tp(\A\tp \Lc)\stackrel{\tau_\A}{\longrightarrow
}\A\tp\A\tp\Lc\tp\Lc
\stackrel{\divideontimes\tp\mathrm{m}}{\longrightarrow }
\A\tp\Lc\tp\Lc \stackrel{\mathrm{m}}{\longrightarrow }
\A\tp \Lc
\nn
\ee
makes $\A\tp \Lc$  an associative $\Ha$-module
algebra, denoted further by $\A\divideontimes\Lc$.
\end{propn}
\begin{proof}
The proof can be conducted by a straightforward verification.
Below we give another proof using our categorical approach to
dynamical associative  algebras, see Example \ref{exDA}.
\end{proof}

\subsection{Infinitesimal analogs of base algebras and dynamical
associative algebras}
\label{ssIA}
In the present subsection, we introduce  quasi-classical analogs
of base algebras and dynamical associative algebras.
\subsubsection{Poisson-Lie manifolds}
Let us recall some basic facts about Poisson-Lie manifolds in

Throughout the text an $\g$-manifold means a manifold equipped
with a left $\g$-action on functions.
This corresponds to a right action on the manifold of a Lie group $G$
relative to $\g$.

Let $\g$ be a Lie bialgebra, i.e. a Lie algebra equipped with a
cobracket map $\mu\colon \g\to \wedge^{\tp 2}\g$ inducing on the
dual space $\g^*$ a Lie algebra structure compatible with the Lie
algebra structure on $\g$ in the sense of  \cite{Dr1}. Recall from
\cite{Dr2} that $\mu$ induces a Poisson structure on the Lie
group $G$ such that the multiplication map $G\times G\to G$ is a
Poisson map (the manifold $G\times G$ is equipped with the
standard Poisson structure of Cartesian product of two Poisson manifolds).
A right
$G$-manifold $P$ is called a Poisson-Lie manifold if the action
$P\times G\to P$ is Poisson. The right $G$-action on $P$ induces a
left action of the universal enveloping algebra $\U(\g)$ on the
function algebra
 $\A(P)$.
For an element $x\in \U(\g)$, let $\vec x_P$ (or simply $\vec x$,
if this causes no confusion) denote the corresponding differential
operator on $P$. For the bidifferential operator on $P$ generated by a
bivector field $\pi$, we use the notation
$\pi(a,b):=(m\circ \pi)(a\tp b)$, $a,b\in \A(P)$,
where $m$ is the multiplication in
$\A(P)$.

The following fact is well known and can be checked directly.
\begin{propn}
Let $\g$ be a Lie bialgebra with cobracket $\mu$, $G$ the
corresponding connected Poisson-Lie group, and $P$ a right
$G$-manifold equipped with a Poisson bracket $\pi$. Then $P$ is a
Poisson-Lie $G$-manifold  if and only if for any $x\in \g$ and
$a,b\in\A(P)$
\be
\label{eqPLB}
\vec x\pi(a,b)-\pi(\vec x a,b)-\pi(a,\vec x b)
=
\overrightarrow{\mu(x)}(a,b).
\ee
\end{propn}
\noindent Any Lie bialgebra structure on $\g$ can be quantized to
a $\C[[t]]$-Hopf algebra $\U_t(\g)$ (quantum group), see \cite{EK}. If
$\A_t(P)$ is a $\U_t(\g)$-equivariant quantization of $\A(P)$,
then quasi-classical limit of $\A_t(P)$ gives a Poisson-Lie bracket
on $P$.

An important particular case of Lie bialgebras is a coboundary one
with the cobracket $\mu(x):=[x\tp 1+ 1\tp x,r]$, where the element
$r\in \wedge^2\g$ satisfies the modified classical Yang-Baxter
equation
\be
\label{cYBE}
[\![r,r]\!]:=[r_{12},r_{13}]+[r_{13},r_{23}]+[r_{12},r_{23}]=\varphi
\in \wedge^3(\g)^\g.
\ee
Formula
(\ref{eqPLB}) then reads
\be
\label{eqPLBr}
[\vec x\tp 1+ 1\tp \vec x, \pi-\vec r]=0.
\ee
In other words, a Poisson-Lie bracket
differs from $\vec r$ by an invariant bivector $f:=\pi-\vec r$ such
that $[\![f,f]\!]$  is equal to $-\vec\varphi$ from (\ref{cYBE}).
Here the operation $f\mapsto[\![f,f]\!]$ is defined by
(\ref{cYBE}) for the Lie algebra of vector fields; this operation is
proportional to the Schouten bracket. Note that the
Poisson-Lie bracket on a Poisson-Lie $\g$-manifold $P$ is the
infinitesimal object for the $\U_t(\g)$-equivariant
quantization of the function algebra on $P$, where $\U_t(\g)$
is the corresponding quantum group. Such brackets are classified
in \cite{DGS1,Kar,D2,DO}
for homogeneous manifolds $G/H$, where $G$ is a simple Lie group
and $H$ its reductive Lie subgroup of maximal rank.

\subsubsection{Poisson base algebras and Poisson base manifolds}
\label{ssPBAPBM}
Let $D(\h)$ denote the double of a Lie bialgebra $\h$,
\cite{Dr1}. As a linear space, $D(\h)$ is the direct sum
$\h+\h^*_{op}$, where $\h^*$ is the dual Lie algebra. The double
$D(\h)$ is endowed with a non-degenerate symmetric bilinear form
induced by the natural pairing between $\h$ and $\h^*_{op}$. There
is a unique extension of the Lie algebra structure from $\h$ and
$\h^*_{op}$ to a Lie algebra on $D(\h)$ such that this
form is ad-invariant. The double is a coboundary Lie bialgebra with the
r-matrix $r:=\sum_i\eta^i\wedge h^i=\frac{1}{2}\sum_i(\eta^i\tp
h^i-h^i\tp \eta^i)$, where $\{h^i\}$ is a basis in $\h$ and
$\{\eta^i\}$ is its dual basis in $\h^*_{op}$. The canonical
element $\theta:=\frac{1}{2}\sum_i(\eta^i\tp h^i+ h^i\tp
\eta^i)$ is ad-invariant. The pair $(r,\theta)$ makes $D(\h)$ a
quasitriangular Lie bialgebra.

\begin{definition}
\label{def_BM} A commutative $D(\h)$-algebra $\Lc_0$ is called
\select{Poisson base algebra} over $\h$, or simply an $\h$-base algebra,
if $\theta$ induces the zero bidifferential operator on $\Lc_0$.

When a Poisson base algebra $\Lc_0$ over $\h$ appears as the function
algebra\footnote{By the function algebra $\A(P)$ we understand,
depending on a particular type of the manifold $P$, the algebra of
polynomial,
analytical, meromorphic, or smooth functions.} on a manifold $L$,
i.e. $\Lc_0:=\A(L)$, we call $L$ a \select{Poisson base manifold} over
$\h$, or simply an $\h$-base manifold.
\end{definition}

\begin{propn}
\label{PbPLm}
An $\h$-base manifold $L$ is a Poisson-Lie $D(\h)$-manifold with
respect to the bracket
\be
\varpi:=\sum_i\vec{\eta^i}\wedge \vec{h^i},
\ee
which is automatically equal to the bivector field $\sum_i \vec{\eta^i}\tp
\vec{h^i}$.
\end{propn}
\begin{proof}
The element $\sum_i{\eta^i}\wedge {h^i}\in \wedge^2 D(\h)$
satisfies the modified Yang-Baxter equation (\ref{cYBE}) with
$\varphi:=[\theta_{12},\theta_{23}]$. Since $\theta$ yields the
zero bivector field on $L$, the three-vector field induced by
$[\theta_{12},\theta_{23}]$ is zero, too. This implies the
following two assertions. Firstly, the bivector $\varpi$ defines a
Poisson structure on $L$. Secondly, any $D(\h)$-invariant Poisson
bracket, hence $\varpi$, is automatically a Poisson-Lie one.
\end{proof}

The following are examples of Poisson base manifolds.
According to Theorem \ref{BAquantization} below,
they can be quantized to $\U_t(\h)$-base algebras, where $\U_t(\h)$ is the
quantized
universal enveloping algebra of $\h$.
\begin{example}[Group spaces $H^*$ and $H$]
\label{H*H}
Let $H$ be  the Lie subgroup in the double $D(H)$ corresponding
to the Lie subalgebra $\h\subset D(\h)$. We will show that the left
coset space $H\backslash D(H)$ is an $\h$-base manifold. Note that
the manifold $H\backslash D(H)$ is naturally identified with the Lie group
$H^*$ corresponding to the Lie algebra $\h^*$. The algebra
$\A(H^*)$ is realized as a subalgebra of functions $f\in
\A\bigl(D(H)\bigr)$ obeying $f(hx)=f(x)$ for $h\in H$. This
subalgebra is invariant under the right regular action of $D(H)$ on
itself. The element $\theta$ is $D(\h)$-invariant, hence the
bivector $\theta^{l,l}-\theta^{r,r}$, where the superscripts $l,r$
denote the left- and right-invariant field extensions, gives the zero
operator on  $\A\bigl(D(H)\bigr)$. Therefore the bivector $\theta^{l,l}$,
which
is equal to $\theta^{l,l}-\theta^{r,r}$ on the left $H$-invariant
functions, gives the zero operator on $\A(H^*)$. Thus, $H^*$ is a
Poisson $\h$-base manifold. In this example, the Poisson bracket on $H^*$
is the
Drinfeld-Sklyanin bracket projected from $D(H)$
to $H^*$.

Similarly to the coset space $H\backslash D(H)$, one can consider the coset
space $H^*\backslash D(H)$, which
is isomorphic as a manifold to the group space $H$. So $H$ is
a Poisson $\h$-base manifold as well.
\end{example}
\begin{example}[Coset spaces $K\backslash H$]
Let us generalize Example \ref{H*H}.
Suppose that $\k$ is a Lie sub-bialgebra in $\h$.
Then the linear sum $\k+\h^*_{op}$ is a Lie sub-bialgebra in $D(\h)$.
Let $K$ be the Lie subgroup in $H$ corresponding to $\k$.
Using the same arguments as in  Example \ref{H*H}, one can
prove that the coset space $K\backslash H$ is a Poisson $\h$-base
manifold.
Indeed, consider the functions on the group $D(H)$ that are invariant
under
the left shifts by elements from $H^*$ and $K$. Such functions
can be identified with functions on $K\backslash H$. The rest
of the construction is exactly the same as in the previous example.
Namely, one can check that the restriction of the Drinfeld-Sklyanin
bracket
makes $K\backslash H$  a Poisson $\h$-base manifold.

It follows that the quotient spaces of the standard Drinfeld-Jimbo simple
Poisson-Lie group $H$ by the Levi and parabolic subgroups
are Poisson $\h$-base manifolds.

Obviously, the same construction holds for the dual Lie bialgebra
$\h^*_{op}$
and its sub-bialgebras; the corresponding coset spaces will be $\h$-base
manifolds.

\end{example}
\begin{example}[Group $H$, the quasitriangular case]
\label{baseH}
Suppose that $\h$ is a quasitriangular Lie
bialgebra, i.e. $\h$ is endowed with an r-matrix $r$ and a
symmetric invariant element $\omega\in \h\tp \h$ such that $r$ satisfies
(\ref{cYBE}) with $\ff:=[\omega_{12},\omega_{23}]$. We can treat
$r$ and $\omega$ as linear maps from $\h^*_{op}$ to $\h$ via
pairing with the first tensor factor. Consider the Lie group $H$
corresponding to $\h$ as a right $H$-manifold via the action
$x\mapsto y^{-1}xy$, $x,y\in H$. This action generates the action
of $\h$ on the function algebra $\A(H)$ by vector fields $\vec
h:=h^l-h^r$, $h\in \h$, where the superscripts $l,r$ stand for the
left- and right- $H$-invariant vector fields generated,
respectively, by the right and the left regular actions of $H$ on
itself. The group $H$  is also a right $\h^*_{op}$-manifold. Namely,
the element $\eta\in \h^*_{op}$
acts on functions from $\A(H)$ as the vector field
$\vec \eta:=r(\eta)^l-r(\eta)^r+\omega(\eta)^l+\omega(\eta)^r$. We have
\be 2 \vec
\theta&=&(r^{l,l}-r^{r,l}-r^{l,r}+r^{r,r})+(\omega^{l,l}-\omega^{r,l}+\omega^{l,r}-\omega^{r,r})
\nn\\
&-&(r^{l,l}-r^{l,r}-r^{r,l}+r^{r,r})+(\omega^{l,l}-\omega^{l,r}+\omega^{r,l}-\omega^{r,r})
\nn\\
             &=&2(\omega^{l,l}-\omega^{r,r}),
\ee which is zero, because $\omega$ is invariant. These actions of
$\h$ and $\h^*_{op}$ define an action of the double
$D(\h)$, thus the group space $H$ is an $\h$-base manifold. In
this example, $\varpi$ is the reflection equation Poisson bracket,
\cite{Sem}. Quantization of this bracket is an RE algebra, cf. Example
\ref{exREa}.

\end{example}
\subsubsection{Poisson dynamical algebras}
In this subsection, we define a Poisson dynamical bracket as an
infinitesimal object for the deformation quantization of a
commutative algebra $\A$ to a dynamical associative algebra,
in the sense of  Definition \ref{defDA}. We assume
that $\A:=\A(P)$, a function algebra on a manifold $P$.

Given a linear space $X$, let $ \mathrm{Alt}$ denote a linear map
from $\End(X^{\tp 3})$ acting by
$$
 \mathrm{Alt}\colon x_1\tp x_2 \tp x_3\mapsto x_1\tp x_2 \tp x_3-x_2\tp x_1
\tp x_3+x_2\tp x_3 \tp x_1,
\quad x_i\in X.
$$
\begin{definition}
\label{def_DB} Let $\h$ be a Lie bialgebra with the cobracket $\mu$,
$L$ an $\h$-base manifold, and $P$ an $\h$-manifold. Let $T(P)$
denote the tangent space to $P$. A function $\pi\colon L \to
\wedge^2 T(P)$ is called a \select{Poisson dynamical bracket} on $P$ (or
on
$\A(P)$) over the base manifold $L$ (or over the Poisson $\h$-base algebra
$\A(L)$) if
\begin{enumerate}
\item
for any $h\in \h$ and $a,b\in \A(P)$
\be
\label{eqDPL}
\vec h_L\pi(\la)(a,b) + \vec h_P\pi(\la)(a,b)-\pi(\la)(\vec h_P
a,b)-\pi(\la)(a,\vec h_P b)= \overrightarrow{\mu(h)}_P(a,b),
\ee
\item
$\pi$ satisfies the equation
\be
\label{eqcDBr}
 \sum_{i} \mathrm{Alt}\bigl(\vec {h^i}_P\tp \vec {\eta^i}_L \pi(\la)\bigr)
 = [\![\pi(\la),\pi(\la)]\!].
\ee
\end{enumerate}
\end{definition}
\noindent
In this definition the expression $\pi(a,b)$ is a
function on $P\times L$. The vector fields $\vec h_L$ and $\vec
h_P$ are induced by the actions of $\h$ on $L$ and $P$,
respectively; $\overrightarrow{\mu(h)}_P$ is a bivector field
induced on $P$ by the element $\mu(h)\in \wedge^2\h$. The
vector field $\vec\eta^i_L$ is induced by the actions of
$\h^*_{op}$ on $L$ (recall that $L$ is a $D(\h)$-manifold).

When $P$ is endowed with a Poisson dynamical bracket over a base
$L$, we say that $\A(P)$ is a \select{Poisson dynamical algebra}.
While a Poisson base manifold is a classical analog of a base algebra,
a Poisson dynamical algebra is a classical analog of dynamical associative
algebra. The Poisson dynamical bracket may be viewed as a map
$\pi\colon \A(P)\wedge \A(P)\to \A(P)\tp \A(L)$.

The following proposition is a generalization of Theorem \ref{Xu}.
\begin{propn}
\label{prop_DBAB0}
Let $\h$ be a Lie bialgebra with cobracket
$\mu$, $L$ an $\h$-base and $P$ an $\h$-manifold. A function
$\pi\colon L \to \wedge^2 T(P)$ is a dynamical bracket on $P$ over
the base $L$ if and only if the bivector
\be
\label{eqLMbr}
\sum_i\vec{\eta^i}_L\wedge \vec {h^i}_L +2\sum_i
\vec{\eta^i}_L\wedge \vec {h^i}_P + \pi
\ee is a Poisson-Lie
bracket on the $\h$-manifold $P\times L$.
\end{propn}
\begin{proof}
This statement is proven by a direct computation. It can be considered as
an infinitesimal analog of Theorem \ref{propCSP}.
\end{proof}

If the base manifold $L$ has $\h$-stable points, then a Poisson dynamical
bracket
on $P$ can be restricted to the coset space $P\slash H$, where it
becomes an ordinary Poisson bracket. This is formalized by the following
proposition.
\begin{propn}
\label{prop_DBAB}
Let $\la_0\in L$ be a stable point under the action of $\h$. Then
$\pi(\la_0)$ induces a Poisson
bracket on the subalgebra of $\h$-invariants
in $\A(P)$.
\end{propn}
\begin{proof}
By the equivariance condition (\ref{eqDPL}),
the function $\pi(\la_0)(f,g)$ is $\h$-invariant when $f,g\in \A(P)$ are
$\h$-invariant.
The Schouten bracket of $\pi(\la_0)$ with itself vanishes on $\h$-invariant
elements from $ \A(P)$,
as follows from (\ref{eqcDBr}).
\end{proof}
\noindent
Proposition \ref{prop_DBAB} gives a quantization method developed
in Section \ref{sQVB} for the
class of Poisson structures coming from Poisson dynamical structures. This
method  uses dynamical associative
algebras, which are quantizations of Poisson dynamical algebras.
\subsection{Quantization of Poisson base algebras and
Poisson dynamical algebras}
Let $\h$ be a Lie bialgebra
and $\U_t(\h)$ the corresponding quantization of $\U(\h)$. Suppose
$\Lc_0$ is a Poisson $\h$-base algebra. By Proposition \ref{PbPLm}, $\Lc_0$
is
endowed with a Poisson bracket induced by the tensor $\sum_i
\eta^i\tp h^i$, where $\{h^i\}\subset \h$ and $\{\eta^i\}\subset
\h^*_{op}$ are dual bases.
\begin{definition}
\label{QofBa} A quantization of the Poisson $\h$-base algebra
$\Lc_0$ is a base algebra $\Lc_t$ over $\U_t(\h)$ that is a
$\U_t(\h)$-equivariant deformation quantization of $\Lc_0$ with
the multiplication
\be
a*_tb=ab+O(t) ,\quad a*_tb-b*_t a= t\sum_i(\vec\eta^i a)(\vec h^i
b)+O(t^2)
\ee
and the coaction $\Lc_t\to\U_t(\h)\ot\Lc_t$ of the form
\be
\delta(a)= 1\tp a + t \sum_i h^i\tp (\vec\eta^i a)+O(t^2),
\label{eqQCcoaction}
\ee
where $a,b\in\Lc_t$.
\end{definition}

When
$\Lc_0=\A(L)$, the function algebra on an $\h$-base manifold, one may
require in the definition that $\Lc_t$ is a star product. Then
$D(\h)$ acts on $\Lc_0$ by vector fields.
\begin{thm}
\label{BAquantization}
Any Poisson base algebra can be quantized.
\end{thm}
\begin{proof}
Let $\theta=\frac{1}{2}\sum_i(h^i\tp \eta^i+\eta^i\tp h^i)\in D(\h)\tp
D(\h)$ be the canonical symmetric invariant of the
double Lie algebra $D(\h)$.
Consider the quasi-Hopf algebra $\U\bigl(D(\h)\bigr)[[t]]$ with the
R-matrix $e^{t\theta}$ and  the  associator $\Phi_t$, which is expressed
through
$t\theta_{12}$ and $t\theta_{23}$, \cite{Dr3}.
Since $\theta$ vanishes on $\Lc_0$, so do $e^{t\theta}$ and $\Phi_t$.
Therefore $\Lc_0[[t]]$ is a commutative algebra not only in the classical
monoidal category of $\U\bigl(D(\h)\bigr)[[t]]$-modules, but also
in the category with the associator $\Phi_t$, i.e. $\Lc_0$ is
$e^{t\theta}$-commutative and  $\Phi_t$-associative.

According to  \cite{EK}, there exists a twist $J_t$ converting the
quasi-Hopf algebra $\U\bigl(D(\h)\bigr)[[t]]$ into
a Hopf one, $\U_t\bigl(D(\h)\bigr)$. This Hopf algebra contains the
quantized enveloping
algebras $\U_t(\h)$ and $\U_t(\h^*_{op})$ as Hopf subalgebras.
The Hopf algebra $\U_t\bigl(D(\h)\bigr)$ is quasitriangular, with the
universal R-matrix
$\Ru_t=(J_t)^{-1}_{21} e^{t\theta}J_t$ lying
in $\U_t(\h^*_{op})\tp \U_t(\h)\subset \U_t\bigl(D(\h)\bigr)\tp
\U_t\bigl(D(\h)\bigr)$.

Applying the twist $J_t$ to the algebra $\Lc_0[[t]]$, we obtain a
quasi-commutative algebra $\Lc_t$ in the
category of $\U_t\bigl(D(\h)\bigr)$-modules. Let us introduce on $\Lc_t$ a
structure of $\U_t(\h)$-comodule
algebra by setting
\be
\label{coaction}
\delta(\ell):=(\Ru_t)_2\tp (\Ru_t)_1\tr \ell, \quad\ell\in \Lc_t.
\ee
Together with the $\U_t(\h)$-action restricted from
$\U_t\bigl(D(\h)\bigr)$,
the coaction (\ref{coaction}) makes $\Lc_t$ a $\U_t(\h)$-base algebra.
This
follows from Corollary \ref{corR-coaction}, where one should set
$\Ha=\U_t\bigl(D(\h)\bigr)$
and $\K=\U_t(\h)$.
\end{proof}
\begin{definition}
\label{def_QCDA}
Let $L$ be an $\h$-base manifold. Let $P$ be an $\h$-manifold and $\pi$ a
Poisson
dynamical bracket on $P$ over $L$.
A quantization of  Poisson dynamical $\h$-algebra
$\A(P)$ is a pair $\bigl(\Lc_t,\A_t(P)\bigr)$, where a)
 $\Lc_t$ is a quantization
of the Poisson base algebra $\A(L)$ in the sense of Definition
\ref{QofBa} and b) $\A_t(P)$ is a flat $\C[[t]]$-module and a dynamical
associative $\U_t(\h)$-algebra
over  $\Lc_t$ such that $\A_t(P)/t\A_t(P)=\A(P)$ and
$a\divideontimes b - b\divideontimes a= t \pi(a,b)+O(t^2)$.
\end{definition}
\begin{conjecture}
\label{APAQ}
Any Poisson dynamical algebra can be quantized.
\end{conjecture}

In Subsection \ref{QVB}, we develop a method of quantizing vector bundles
on  the coset space $P/H$,
using dynamical associative algebras.
By duality, the construction of Section \ref{QVB}
can be formulated in terms of base algebras rather than coalgebras.
Namely, let $\A$ be a dynamical associative algebra over an $\Ha$-base
algebra $\Lc$ and let $\chi$ be an $\Ha$-invariant character of $\Lc$.
Then the composition map
$\A\tp \A\stackrel{\divideontimes}{\longrightarrow} \A\tp \Lc \stackrel{\id
\tp \chi}{\longrightarrow} \A$
yields an associative multiplication on the subspace of $\Ha$-invariant
elements of $\A$.
Thus invariant characters of base algebras are important for our approach
to quantization,
see also \cite{DM1}.

In the deformation situation, the infinitesimal analogs of
$\U_t(\h)$-invariant characters of the base algebra $\Lc_t$ are
$\h$-stable points on the $\h$-base manifold $L$.
By Proposition \ref{prop_DBAB}, each $\h$-stable point defines a Poisson
structure on $P/H$.
It is that Poisson structure, which is quantized by applying an invariant
character to
the dynamical associative quantization of $\A(P)$.
The question is whether every $\h$-stable point can be quantized to a
$\U_t(\h)$-invariant character of $\Lc_t$.
The answer to this question is affirmative.
\begin{propn}
Let $\Lc_t$ be the quantization of the function
algebra on a base manifold $L$ built in Theorem \ref{BAquantization}.
Then every $\h$-stable point $\la_0$ on $L$ defines
an $\U_t(\h)$-invariant character of $\Lc_t$ by
 $\chi^{\la_0}(f)= f(\la_0)$ for $f\in \Lc_t$.
\end{propn}
\begin{proof}
As follows from the explicit form of the twist $J_t$ built in \cite{EK},
it reduces to $1\tp 1$ at every $\h$-stable point $\la_0\in L$.
It follows from the proof of  Theorem \ref{BAquantization} that
the star product in $\Lc_t$ satisfies the condition
$(f*g)(\la_0)=(fg)(\la_0)=f(\la_0) g(\la_0)$
for any pair of functions $f,g\in \A(L)$.
\end{proof}

\subsection{Dynamical Yang-Baxter equations}
\label{ssDYBEs}
In this subsection we give  definitions of the
classical and quantum dynamical Yang-Baxter equations over an arbitrary
base algebra.

\begin{definition}
\label{defDrm}
Let $\g$ be a Lie bialgebra and $\h\subset \g$ its sub-bialgebra; let $\mu$
denote
the Lie cobracket on $\h$.
Let $L$ be a Poisson $\h$-base manifold. A function $\bar
r\colon L\to \g\tp\g$ is called a \select{classical dynamical r-matrix}
over base $L$ if
\begin{enumerate}
\item
for any $h\in\h$
\be
\label{inv}
\vec h_L \bar r(\la) + [h\tp 1+ 1\tp h,\bar r(\la)] = \mu(h),
\ee
\item
the sum $\bar r(\la)+\bar r_{21}(\la)$ is $\g$-invariant,
\item
$\bar r$ satisfies the  equation
\be \label{cDYBE}
\sum_{i}
\mathrm{Alt}\bigl(h^i\tp \vec {\eta^i}_L \bar r(\la)\bigr)
 = [\![\bar r(\la),\bar r(\la)]\!].
\ee
\end{enumerate}
\end{definition}
\noindent We call (\ref{cDYBE}) the classical dynamical
Yang-Baxter equation over the base $L$.
Condition (\ref{inv}) means  quasi-equivariance of $\bar r(\la)$
with respect to the action of $\h$.
In fact, the symmetric part
$\bar\omega=\frac{1}{2}(\bar r+\bar r_{21})$
is constant on every $D(\h)$-orbit in $L$, i.e. $\vec x_L\bar \omega=0$
for any $x\in D(\h)$.

Consider the opposite Lie bialgebra $\h_{op}$ equipped with the opposite
bracket and
the same cobracket.
Equipped with the   opposite bracket, the manifold $L$ becomes a Poisson
base manifold for $\h_{op}$.
The Cartesian product $L\times L$ is a Poisson base manifold with
respect to the Lie bialgebra $\h\oplus \h_{op}$.
Let $G$ be the Lie group corresponding to $\g$.
The following proposition characterizes the dynamical r-matrices.
\begin{propn}
\label{G-PDA0}
A function $\bar r\colon L\to \g\tp \g$ is a dynamical r-matrix over the
Poisson $\h$-base manifold $L$
if and only if $M:=G\times L\times L$ is equipped with the $\h\oplus
\h_{op}$- Poison Lie structure
such that  the projection $M\to L\times L$ is a Poisson map and
\be
\{f,a\}&:=&\sum_i(\vec\eta^i_L f)\bigl(h^{i\:l} + h^{i\:r}\bigr)(a)
\\
\{a,b\}&:=&\bigl(\bar r^{l,l}(\la')-\bar r^{r,r}(\la'')\bigr)(a,b)
\ee
for $f\in \A(L\times L)$, $a,b\in \A(G)$, and $(\la',\la'')\in L\times L$.
\end{propn}
\begin{proof}
Straightforward.
\end{proof}

Suppose that $\g$ is a quasitriangular Lie bialgebra with  an r-matrix
$r_\g\in \g\tp \g$.
Let $\omega=\frac{1}{2}(r_\g+r^{21}_\g)$ denote the
symmetric part of $r_\g$.
Assume that $L$ is $D(\h)$-transitive, i.e. $D(\h)$-invariants in $\A(L)$
are scalars.
\begin{propn}
\label{G-PDA}
A function $\bar r\colon L\to \g\tp \g$ subject to $\frac{1}{2}(\bar r+\bar
r)=\omega\in \g\tp \g$
is a dynamical r-matrix
if and only if $M:=L\times G$ is equipped with
an $H$-invariant Poisson structure, such that the projection
$M\to L$ is a Poisson map, and
\be
\{f,a\}:=\sum_i(\vec\eta^i_L f)\bigl(h^{i\:l} a\bigr)
,\quad
\{a,b\}:=(\bar r^{l,l}-r^{r,r}_\g)(a,b)
\ee
for $f\in \A(L)$, $a,b\in \A(G)$,
\end{propn}
\begin{proof}
Straightforward.
\end{proof}
Proposition \ref{G-PDA0} implies that  the bivector field $\bar
r^{l,l}(\la')-\bar r^{r,r}(\la'')$ makes $\A(G)$ a
Poisson dynamical algebra over the $\h\oplus \h_{op}$-base manifold
$L\times L$. By Proposition \ref{G-PDA},
the bivector field $\bar r^{l,l}(\la)-r^{r,r}_\g$ makes $\A(G)$ a
Poisson dynamical algebra over the $\h$-base manifold $L$.

\begin{propn}
Let $\bar r\colon L\to \g\tp \g$ be a classical dynamical r-matrix on an
$\h$-base manifold $L$.
Let $r_\g\in \g\tp \g$ be a constant r-matrix whose symmetric part
coincides with the symmetric part
of $\bar r$.
Suppose that  $\la_0\in L$ is  an $\h$-stable point.
Then the bivector field $\bar r^{l,l}(\la_0)-r^{r,r}_\g$ yields a
$\g$-Poisson-Lie structure
on the coset space $G\slash H$.
\end{propn}
\begin{proof}
Applying Proposition \ref{prop_DBAB} to the Poisson dynamical bracket $\bar
r^{l,l}(\la)-r^{r,r}_\g$
on $G$, we obtain a Poisson structure on the subalgebra in $\A(G)$ that
consists of invariants under the action of $\h$ by the left-invariant
vector fields.
This algebra is canonically identified with the algebra of functions on the
coset
space $G\slash H$. Obviously, this Poisson structure is a Poisson-Lie one,
with respect to the $\g$-action induced by the left $G$-action.
\end{proof}
\begin{remark}
Let $\h=\g$ be quasitriangular, with the classical r-matrix $r_\g$ whose
symmetric part is equal to that of $\bar r$.
Then $\bar r(\la)-r_\g$ is the dynamical r-matrix of \cite{FhMrsh}.
\end{remark}

We complete this subsection with a definition of the
quantum dynamical Yang-Baxter equation over an arbitrary base algebra.
This definition naturally follows from our categorical point of view
presented in Section \ref{sDCRA}.

The quantum DYBE will be defined for any triple $(\U,\Ha,\Lc)$,
where $\Ha$ is a Hopf subalgebra in a Hopf algebra $\U$ and  $\Lc$ is an
$\Ha$-base algebra,
\begin{definition}
\label{def_QDR} An element $\bar
\Ru=\bar\Ru_1\tp\bar\Ru_2\tp\bar\Ru_3 \in \U\tp \U\tp \Lc$ is
called a \select{universal quantum dynamical R-matrix} of $\U$ over the
$\Ha$-base algebra $\Lc$ if it satisfies the equivariance condition
\be
h^{(2)}\bar\Ru_1\tp h^{(1)}\bar\Ru_2\tp h^{(3)}\tr \bar\Ru_3
&=&
\bar\Ru_1 h^{(1)}\tp \bar\Ru_2 h^{(2)}\tp \bar\Ru_3, \quad h\in \Ha,
\label{uRinv}
\ee
and the quantum dynamical Yang-Baxter
equation
\be \bar \Ru_{12} \;{}^{(2)}\!\bar \Ru_{13}\; \bar
\Ru_{23}&=&
 \;{}^{(1)}\!\bar  \Ru_{23} \;\;\bar  \Ru_{13}  \;{}^{(3)}\!\bar \Ru_{12},
\label{uR}
\ee
in $\U\tp\U\tp \U\tp \Lc$.
\end{definition}
\noindent
Here  the notation ${}^{(i)}\!\bar  \Ru$ means the following.
Applying the coaction to the $\Lc$-component of $\Ru$,
we get the element
${}^{(3)}\!\bar  \Ru:= \bar\Ru_1\tp \bar\Ru_2 \tp \bar\Ru_3^{(1)}\tp
\bar\Ru_3^{[2]}$.
The other two are obtained from this by permutations, namely
${}^{(2)}\!\bar  \Ru:= \bar\Ru_1\tp \bar\Ru_3^{(1)}\tp \bar\Ru_2 \tp
\bar\Ru_3^{[2]}$
and
${}^{(1)}\!\bar  \Ru:= \bar\Ru_3^{(1)}\tp \bar\Ru_1\tp \bar\Ru_2 \tp
\bar\Ru_3^{[2]}$.

Equation (\ref{uR}) specializes to (\ref{XuYBE}) for
$\Ha=\U(\h)$, $\U=\U(\g)$, and $\Lc$ being the extension of
$\U(\h_t)$ to the PBW star product on functions on $\h^*$, cf. Example
\ref{PBW*pr}.
In this case, equation  (\ref{uR}) coincides with the conventional
dynamical
Yang-Baxter equation (\ref{DYBE})  for $\h$ a commutative Lie
subalgebra in $\g$ and $\Lc$ being the algebra of functions
on $\h^*$, \cite{EV2}.

Suppose that $\U$ and $\Ha$ are quantizations of the universal enveloping
algebras
$\U(\g)$ and $\U(\h)$ and $\Lc$ is a quantization of the function algebra
on a
Poisson $\h$-base manifold $L$.
Suppose that the universal dynamical R-matrix has the form
$\Ru=1\tp 1 \tp 1 + t \bar r + O(t^2)$.
Then $\bar r$ is a function on $L$ with values in $\g\tp \g$. It satisfies
the equations (\ref{inv}) and (\ref{cDYBE}), which are the consequences
of the equations (\ref{uRinv}) and (\ref{uR}).

\begin{remark}
The definitions of classical and quantum dynamical R-matrix given above
admit further generalization.
The reader is referred to \cite{DM5}, where the classical dynamical
r-matrices
are studied in connection with Lie bialgebroids. The present definitions
are
conditioned by our specific approach confined to
the strict monoidal categories (i.e. with trivial associator). If one
considers general monoidal
categories, as in Example \ref{Phi-twistor}, the equation (\ref{uR}) would
involve an associator. In
the quasi-classical limit, it gives the dynamical r-matrix of
the Alekseev-Meinrenken type, \cite{AM}.
\end{remark}
\section{Dynamical categories}
\label{sDC}
\subsection{Base algebra in a monoidal category}
A  dynamical associative  algebra $\A$ from Definition
\ref{defDA} may serve as a model for further generalizations. It
turns out that there is a monoidal category where $\A$ is an associative
algebra. Such categories can be built for all Hopf algebras and
they include dynamical categories of Etingof-Varchenko introduced for
commutative cocommutative
Hopf algebras in \cite{EV3}. Such notions
as dynamical twist and dynamical Yang-Baxter equation can be
naturally formulated and generalized within the dynamical
categories, which are subjects of our further study.

Let $\hat\O$ be a monoidal category.
We will work, for simplicity, with only strict monoidal categories, i.e.
having the trivial associator;  all the constructions can be carried over
to the general
case in a straightforward way. Let $Z(\hat\O)$ be the center of  $\hat\O$,
see \cite{Kas}.
The center is a braided monoidal category consisting of pairs
$(A,\tau)$, where $A$ is an object of $\hat \O$ and $\tau$
the collection of permutations $\tau_X\colon A\tp X\to X\tp A$, satisfying
natural conditions.
\begin{definition}
A \select{base algebra} in the category $\hat\O$ is commutative algebra
from
a $Z(\hat \O)$.
\end{definition}

In other words, a  base algebra is an algebra in $\hat \O$ and a collection
of morphisms
$\tau_A\in \Hom_{\hat\O}(\Lc\tp A,A\tp \Lc)$, $A\in \Ob\; \hat\O$,
such that the  following diagrams are commutative:
\be
&
   \dgARROWLENGTH=0.8\dgARROWLENGTH
\begin{diagram}
\node{\Lc\tp
B}\arrow{e,t}{\tau_{B}}\arrow{s,l}{\id_\Lc\tp\psi}\node{B\tp\Lc}\arrow{s,r}{\psi\tp\id_\Lc}
\\
\node{\Lc\tp A}\arrow{e,t}{\tau_{A}}\node{A\tp\Lc}
\end{diagram}
\label{si_nat}
\\[8pt]
&
\begin{diagram}
   \dgARROWLENGTH=0.1\dgARROWLENGTH
\node{\Lc\tp A\tp B}\arrow[2]{e,t}{\tau_{A\tp
B}}\arrow{se,b}{\tau_{A}}\node[2]{A\tp B\tp\Lc}\arrow{sw,r}{\tau_{B}}
\\
\node[2]{A\tp \Lc\tp B}
\end{diagram}
\label{si_hom}
\\[8pt]
&
\begin{diagram}
   \dgARROWLENGTH=0.6\dgARROWLENGTH
\node{\Lc\tp \Lc\tp
A}\arrow{e,t}{\tau_{A}}\arrow{s,l}{\mathrm{m}_\Lc}\node{\Lc\tp A\tp
\Lc}\arrow{e,t}{\tau_{A}}\node{A\tp\Lc\tp\Lc}
\arrow{s,r}{\mathrm{m}_\Lc}
\\
\node{\Lc\tp A}\arrow[2]{e,t}{\tau_{A}}\node[2]{A\tp\Lc}
\end{diagram}
\label{si_hom1}
\\[8pt]
&
\begin{diagram}
   \dgARROWLENGTH=0.5\dgARROWLENGTH
\node{\Lc\tp
\Lc}\arrow[2]{e,t}{\tau_{\Lc}}\arrow{se,b}{\mathrm{m}_\Lc}\node[2]{\Lc\tp\Lc}\arrow{sw,r}{\mathrm{m}_\Lc}
\\
\node[2]{\Lc}
\end{diagram}
\label{si_L}
\ee
for
all $A,B\in \Ob\;\hat\O$, $\psi\in \Hom_{\hat\O}(B,A)$.
\begin{example}
The unit object $1_{\hat \O}$ is the simplest example of a base algebra.
The algebra structure and permutation are defined by the canonical
isomorphisms
$1_{\hat \O}\tp A\simeq A\simeq A\tp 1_{\hat \O}$ for all $A\in \Ob
\;\hat\O$.
\end{example}

\begin{example}
\label{qcomalg}
When the category $\hat\O$ is braided with braiding $\si$, any commutative
algebra $\Lc$
in this category has two natural base algebra structures, with respect to
the $\tau=\si$ and $\tau=\si^{-1}$.
\end{example}

\begin{example}
\label{exO_H}
Let $\Ha$ be a Hopf algebra and $\hat\O$ the monoidal category of left
$\Ha$-modules.
Any base $\Ha$-algebra in the sense of Definition \ref{def_BA} is a base
algebra in the category $\hat \O$.
Indeed, for a left $\Ha$-module, $A$, we define the permutation
$\tau_A\colon  \Lc\tp A\to A\tp \Lc$
by
\be
\label{si_inv}
\ell\tp a\to \ell^{(1)}\tr a \tp \ell^{[2]}, \quad a\in A, \> \ell\in \Lc.
\ee
The permutation (\ref{si_inv}) is $\Ha$-equivariant, as follows from
(\ref{eqDYM}), hence the condition
(\ref{si_nat}) is satisfied.
Conditions (\ref{si_hom}) and (\ref{si_hom1}) hold because $\Lc$ is an
$\Ha$-comodule algebra.
Equation (\ref{si_L}) follows from (\ref{eqQCDYA}).
\end{example}

\begin{example}
Let $\hat\O$ be the category of semisimple modules
over a commutative finite dimensional Lie algebra $\h$.
Take $\Lc$ to be the algebra of
functions
on $\h^*$, which is a trivial $\h$-module. Let $A$ be an semisimple
$\h$-module.
The permutation $\tau_A$ between $\Lc$ and $A$
is defined by $f(x)\tp a\mapsto a \tp f\bigl(x+\al(a)\bigr)$,
where $f\in\Lc$, $a\in A$, and $\al(a)$ is weight of $a$.
\end{example}

\subsection{Dynamical categories over base algebras}
\label{ssExtL}
Let $(\hat \O,\tp)$ be a monoidal category and $\O$ be a monoidal
subcategory in $\hat \O$.
Given a base algebra $(\Lc,\tau)$ in $\hat \O$, let us
construct a new monoidal category
$\bar \O_{\Lc}$. Objects in
$\bar \O_{\Lc}$ are the same as in  $\O$.
For two objects $A$ and $B$ in $\bar \O_{\Lc}$, morphisms $\Hom_{\bar
\O_{\Lc}}(A,B)$
are $\hat\O$-morphisms $\Hom_{\hat\O}(A, B\tp \Lc)$.
Since the algebra $\Lc$ is unital, every morphism $\phi\in \Hom_{\O}(A,B)$
naturally becomes a morphism from $\Hom_{\bar \O_{\Lc}}(A,B)$
through the composition
$A\stackrel{\phi}{\longrightarrow}{B\tp 1_\O}\to  {B\tp \Lc}$.

The composition of two morphisms $A\stackrel{\phi}{\longrightarrow} B$ and
$B\stackrel{\psi}{\longrightarrow} C$
in $\bar \O_{\Lc}$ is defined as the composition
\be
A\stackrel{\phi}{\longrightarrow} B\tp \Lc \stackrel{\psi}{\longrightarrow}
  C\tp \Lc\tp \Lc
\stackrel{\mathrm{m}}{\longrightarrow} C\tp \Lc,
\ee
in $\hat\O$, where the rightmost arrow is the multiplication in $\Lc$.
It is easy to see that the composition  is associative.
The identity morphism $\id_A$ for $A\in \Ob \; \O_\Lc$ is the
composition $A\to A\tp 1_{\hat \O} \to A\tp \Lc $, where the first arrow
is
the canonical isomorphism and the second one is the natural inclusion
$1_{\hat \O}\to \Lc$ via the unit
of $\Lc$.
Thus  $\bar \O_{\Lc}$ is a category.

Let us introduce a monoidal structure $\bar \tp$ in $\bar \O_{\Lc}$ setting
it on objects
as in  $\O$; on the morphisms it is defined by the composition
\be
\label{eqTP}
A\tp C \stackrel{\phi\tp \psi}{\longrightarrow}
B\tp \Lc \tp D\tp \Lc \stackrel{\tau_D}{\longrightarrow}   B \tp D\tp
\Lc\tp \Lc
\stackrel{\mathrm{m}_\Lc}{\longrightarrow} B \tp D\tp \Lc,
\ee
for $\phi\in\Hom_{\bar \O_{\Lc}}(A,B)$ and $\psi\in\Hom_{\bar
\O_{\Lc}}(C,D)$.
\begin{propn}
The tensor product $\bar \tp$ defined by (\ref{eqTP}) makes $\bar \O_\Lc$
a monoidal category.
\end{propn}
\begin{proof}
The unit object from $\O$ is obviously the neutral element for $\bar\tp$.
Let us prove associativity of $\bar \tp$.
Using compatibility (\ref{si_hom1}) of $\tau$ with the multiplication
$\mathrm{m}_\Lc$
and associativity (\ref{si_hom}) we find that the diagram
\be
   \dgARROWLENGTH=0.7\dgARROWLENGTH
\begin{diagram}
\node{\Lc\tp A\tp\Lc\tp
B}\arrow{s,l}{\tau_B}\arrow{e,t}{\tau_A}\node{A\tp\Lc\tp \Lc\tp
B}\arrow{e,t}{\mathrm{m}_\Lc}\node{A\tp\Lc\tp B} \arrow{s,r}{\tau_B}
\\
\node{\Lc\tp A\tp B\tp\Lc}\arrow{e,t}{\tau_{A\tp B}}\node{A\tp B \tp
\Lc\tp\Lc} \arrow{e,t}{\mathrm{m}_\Lc} \node{A\tp B \tp \Lc}
\end{diagram}
\nonumber
\ee
is commutative for all $A,B\in \Ob\; \O$. From this one can deduce
associativity of $\bar \tp$.

Now we will prove functoriality. It is equivalent
to the four conditions:
\be
(\id \bar \tp \phi)\bar \circ  (\id \bar \tp \psi)&=&\id \bar \tp (\phi\bar
\circ \psi),
\label{eqN2}\\
(\phi\bar \tp \id)\bar \circ (\id \bar \tp \psi)&=&\phi \bar \tp \psi,
\label{eqN4}\\
(\phi\bar \tp \id)\bar \circ  (\psi\bar \tp \id)&=&(\phi\bar \circ
\psi)\bar \tp \id,
\label{eqN1}\\
(\id \bar \tp \psi)\bar \circ  (\phi\bar \tp \id)&=&\phi \bar \tp \psi
\label{eqN3}
\ee
for all pairs of morphisms $\phi$, $\psi$.
Observe that $\phi\bar \tp \id_B=(\id_{A'} \tp \tau_B)\circ(\phi \tp
\id_B)$ and  $\id_B \bar \tp \phi=\id_B \tp \phi$
for any morphism $A\stackrel{\phi}{\longrightarrow}A'$ and any object $B$.
This immediately leads to
(\ref{eqN2}) and (\ref{eqN4}). Condition  (\ref{eqN1}) follows from
(\ref{si_hom1}).
Let us prove condition
(\ref{eqN3})  assuming $\phi\in \Hom_{\bar \O_\Lc}(A,A')$ and $\psi\in
\Hom_{\bar \O_\Lc}(B,B')$.
It suffices to show that the following diagram is commutative (the identity
maps are suppressed):
\be
   \dgARROWLENGTH=0.1\dgARROWLENGTH
\begin{diagram}
\node{A\tp B}\arrow{e,t}{\phi}\arrow{se,b}{\phi\tp \psi}\node{A'\tp\Lc\tp
B}
\arrow{s,r}{\psi}\arrow[2]{e,t}{\tau_B}\node[2]{A'\tp
B\tp\Lc}\arrow{s,r}{\psi}
\\
\node[2]{A'\tp \Lc\tp B'\tp\Lc}\arrow{se,b}{\tau_{ B'}}
\arrow[2]{e,t}{\tau_{B'\tp\Lc}}\node[2]{A'\tp B'\tp\Lc\tp\Lc}
\arrow[2]{s,r}{\mathrm{m}_\Lc}
\\
\node[3]{A'\tp
B'\tp\Lc\tp\Lc}\arrow{ne,t}{\tau_{\Lc}}\arrow{se,b}{\mathrm{m}_\Lc}
\\
\node[4]{A'\tp \B'\tp\Lc}
\end{diagram}
\label{si_qC}
\ee
Commutativity of the rectangle follows from (\ref{si_nat}); the two lower
triangles are commutative
by virtue of (\ref{si_hom}
) and (\ref{si_L}).
\end{proof}
The category $\bar\O_\Lc$ naturally includes $\O$ as a monoidal
subcategory.
We call $\bar \O_\Lc$ the \select{dynamical extension} of $\O$ over the
base algebra $\Lc$.

\begin{example}
The simplest example is when $\Lc=1_{\hat \O}$ and $\O=\hat\O$; then the
category $\bar \O_{\Lc}$ is canonically isomorphic to $\O$.
\end{example}

\begin{example}
\label{exDHL}
Let $\Ha$ be a Hopf algebra and $\hat \O$ the category of left
$\Ha$-modules.
As was mentioned in Example \ref{exO_H}, any $\Ha$-base algebra, including
$\Ha$ itself,  is a base algebra in $\hat \O$.
Let $\Mc_\Ha$ be the subcategory of locally finite
$\Ha$-modules (a module is called locally finite if every its element lies
in
a finite dimensional submodule). Its dynamical extension over a base
algebra $\Lc$ is denoted further
by $\bar\Mc_{\Ha;\Lc}$, or simply $\bar\Mc_{\Ha}$ for $\Lc=\Ha$.
\end{example}
\subsection{Morphisms of base algebras}
By a morphism of base algebras $(\Lc_1,\tau^1)\to (\Lc_2,\tau^2)$ in a
category $\hat\O$
we mean a morphism of
$\hat \O$-algebras
$\Lc_1\stackrel{f}{\longrightarrow}\Lc_2$ such that the diagram
\be
&
   \dgARROWLENGTH=0.8\dgARROWLENGTH
\begin{diagram}
\node{\Lc_1 \tp A}\arrow{e,t}{\tau^1_A}\arrow{s,l}{f\tp
\id_A}\node{A\tp\Lc_1}\arrow{s,r}{\id_A\tp f}
\\
\node{\Lc_2 \tp A}\arrow{e,t}{\tau^2_A}\node{A\tp\Lc_2}
\end{diagram}
\nonumber
\ee
is commutative for all $A\in \Ob \; \hat\O$.
\begin{example}
Let $\Ha$ be a Hopf algebra and $\hat \O$ the category of left
$\Ha$-modules.
A homomorphism of two $\Ha$-base algebras can be defined
as a homomorphism of $\Ha$-algebras and $\Ha$-comodules. Then it
a morphism of base algebras in  $\hat \O$, cf. Example \ref{exO_H}.
\end{example}

\begin{example}
Any invariant character $\chi$  of $\Lc$
defines a homomorphism of  base algebras $\Lc\to \Ha$ by the formula
$\ell\mapsto \ell^{(1)}\chi\bigl(\ell^{[2]}\bigr)$.
Indeed, this is an algebra and coalgebra map because  $\Lc$ is an
$\Ha$-comodule algebra.
This map is equivariant for invariant $\chi$, by virtue of (\ref{eqDYM}).
\end{example}

Recall that a functor from one monoidal category to another
is called strong monoidal if it is unital (relates the units) and
commutes with tensor products.
We conclude this subsection with an obvious proposition.
\begin{propn}
A morphism of base algebras $(\Lc_1,\tau^1)\to (\Lc_2,\tau^2)$ induces a
strong monoidal functor
$\bar \O_{\Lc_1}\to \bar \O_{\Lc_2}$.
\end{propn}
\subsection{Category ${\bar{\mathcal{M}}^{\Ha^*}}$}
\label{ssCO_H}
The dynamical extension of a monoidal category can be defined using
a notion of base coalgebra instead of base algebra.
We will present such a formulation for the case when the monoidal
category $\hat\O$ is a category of $\Ha$-modules and the base coalgebra is
a restricted dual
to $\Ha$.

Let $\Ha^*$ denote the Hopf algebra formed
by matrix elements of finite dimensional  semisimple representations of
$\Ha$
(we assume that the supply
of such elements is big enough to induce a non-degenerate pairing between
$\Ha^*$ and $\Ha$).
We equip $\Ha^*$ with the structure of a left $\Ha$-module with respect to
the action
\be
\label{eqModH*}
x\tp \la \mapsto x^{(2)}\tr \la \tl\gm(x^{(1)}),
\quad x\in \Ha, \quad\la \in \Ha^*,
\ee
expressed through the coregular left and right actions, $\tr$ and $\tl$, of
$\Ha$ on $\Ha^*$.

Let $\hat O$ be the category of left $\Ha$-modules. We can
consider the  category of locally finite
right $\Ha^*$-comodules
as a subcategory in  $\hat\O$, since every right $\Ha^*$-comodule is
naturally a  left $\Ha$-module. We denote this category by $\Mc^{\Ha^*}$.

The following statement introduces a permutation between
$\Ha^*$ and other $\Ha^*$-comodules.
\begin{propn}
\label{sigma}
For any $A\in \Ob\; \Mc^{\Ha^*}$ the map $\tau^{A}\colon \Ha^*\tp{A}\to
{A}\tp \Ha^*$
defined as
\be
\label{tau*}
\tau^{A}(\la\tp a):=a^{[0]}\tp \la a^{(1)}
\ee
is an
isomorphism of $\Ha$-modules.
\end{propn}
\begin{proof}
First of all observe that $\tau^A$ is invertible and its inverse is
$$
(\tau^A)^{-1} (a\tp\la) = \la\gamma^{-1}(a^{(1)})\tp a^{[0]},
\quad \la\in \Ha^*, \> a\in A.
$$
Further, for all $x,y\in \Ha$ we have
\be
\langle \tau^{A}\bigl(x\tr(\la\tp a)\bigr),\id\tp y\rangle
&=&
\langle \tau^{A}\bigl(x^{(1)}\tr\la\tp a^{[0]} \bigr),\id\tp y\rangle
\langle a^{(1)},x^{(2)}\rangle
\hspace{3.4cm}
\nn\\&=&
a^{[0]}\tp \langle (x^{(1)}\tr\la)a^{(1)}, y \rangle \langle
a^{(2)},x^{(2)}\rangle
 \nn\\&=&
a^{[0]}  \tp \langle x^{(1)}\tr \la,y^{(1)}\rangle\langle
a^{(1)},y^{(2)}x^{(2)}\rangle,
\label{eq_00}
\ee
On the other hand,
\be
\langle x\tr\tau^{A}(\la\tp a),\id\tp y\rangle &=&
a^{[0]} \tp \langle a^{(1)},x^{(1)}\rangle\langle \la a^{(2)},\gm
(x^{(2)})yx^{(3)}\rangle
 \nn\\&=&
a^{[0]}  \tp \langle a^{(1)},x^{(1)}\rangle
\langle \la,\gm (x^{(3)})y^{(1)}x^{(4)}\rangle
\langle a^{(2)},\gm(x^{(2)})y^{(2)}x^{(5)}\rangle
\nn\\&=&
a^{[0]}  \tp \langle  a^{(1)},y^{(2)}x^{(3)}\rangle
\langle \la,\gm(x^{(1)})y^{(1)}x^{(2)}\rangle
\label{eq_01}
\ee
for all $x,y\in \Ha$, $\la\in \Ha^*$, $a\in A$.
The resulting expression in (\ref{eq_01}) is easily brought to
(\ref{eq_00}).
\end{proof}

Let us define the dynamical extension, $\bar \Mc^{\Ha^*}$,
of the category $\Mc^{\Ha^*}$.
The objects in $\bar \Mc^{\Ha^*}$ are right $\Ha^*$- comodules.
The set of morphisms $\Hom_{\bar\Mc^{\Ha^*}}(A,B)$
consists of $\Ha$-equivariant maps from ${\Ha^*}\tp A$ to $B$.
The composition $\phi\bar \circ \psi$ of morphisms $\phi\in \Hom(A,A')$ and
$\phi\in \Hom(A',A'')$
is defined as the composition map
\be
\label{eq_compH*}
\begin{diagram}
\node{{\Ha^*} \tp A}\arrow{e,t}{\Delta\tp \id_A}\node{{\Ha^*}\tp {\Ha^*}
\tp A}
\arrow{e,t}{\id_{\Ha^*} \tp \psi}\node{{\Ha^*}\tp
A'}\arrow{e,t}{\phi}\node{A''}
\end{diagram}
.
\ee
This operation is apparently associative and $\ve\tp \id_A$ is
the identity in $\Hom_{\bar\Mc^{\Ha^*}}(A,A)$; here $\ve$ is the counit in
$\Ha^*$.

Now we introduce a monoidal structure on $\bar\Mc^{\Ha^*}$.
We put the tensor product of objects from  $\bar\Mc^{\Ha^*}$
as in  $\Mc^{\Ha^*}$.
The tensor product $\phi\bar\tp \psi$ of $\phi\in
\Hom_{\bar\Mc^{\Ha^*}}(A,A')$ and $\psi\in \Hom_{\bar\Mc^{\Ha^*}}(B,B')$
is defined as the composition
\be
\label{eq_tpH*}
{\Ha^*}\tp A\tp B \stackrel{\Delta}{\longrightarrow}{\Ha^*}\tp{\Ha^*}\tp
A\tp B
\stackrel{\tau^A}{\longrightarrow} {\Ha^*}\tp A\tp {\Ha^*}\tp B
\stackrel{\phi\tp \psi}{\longrightarrow} A'\tp B'.
\ee
One can check that, indeed,  the operation $\bar \tp$ makes
$\bar\Mc^{\Ha^*}$ a monodial category.
\subsection{Comparison of categories $\bar {\Mc}^{\Ha^*}$ and $\bar
{\Mc}^{\Ha^*}_{\Ha}$}
Since $\Mc^{\Ha^*}$ is a subcategory in the category of $\Ha$-modules,
it can be extended
to the dynamical category $\bar\Mc^{\Ha^*}_\Ha$
over the base algebra $\Lc=\Ha$ along the line of Subsection \ref{ssExtL}.
Our next goal is to compare the categories $\bar\Mc^{\Ha^*}_\Ha$
and  $\bar\Mc^{\Ha^*}$. Since they have the same supply of objects,
we will study relations between their morphisms.

Introduce a pairing between $\Ha^*$ and $\Ha$ by the formula
\be
\label{eq_pairing}
(h,x):=\langle \gm^{-1}(h), x\rangle,
\ee
where $\langle.,.\rangle$ is the canonical Hopf pairing.
It is invariant under the adjoint action of $\Ha$ on itself (\ref{eq_ad})
and on $\Ha^*$ given by (\ref{eqModH*}).
\begin{lemma}
\label{lm_t-t}
For any right $\Ha^*$-comodule $A\in \Mc^{\Ha^*}$ the diagram
\be
&
   \dgARROWLENGTH=0.7\dgARROWLENGTH
\begin{diagram}
\node{\Ha^*\tp A \tp \Ha
}\arrow{e,t}{\tau^A}\arrow{s,l}{\tau^{-1}_A}\node{A \tp \Ha^*\tp
\Ha}\arrow{s,r}{(.,.)}
\\
\node{\Ha^*\tp \Ha \tp A}\arrow{e,t}{(.,.)}\node{A}
\end{diagram}
\ee
is commutative.
\end{lemma}
\begin{proof}
Straightforward.
\end{proof}
To any equivariant map $\phi\colon A\to B\tp \Ha$ we put into
correspondence
an equivariant map $\phi'\colon\Ha^*\tp A\to B$ being the composition
\be
\label{eq_HH*}
\Ha^*\tp A\stackrel{\phi}{\longrightarrow} \Ha^*\tp B\tp \Ha
\stackrel{\tau^B}{\longrightarrow}
 B\tp \Ha^*\tp\Ha \stackrel{(.,.)}{\longrightarrow}B.
\ee
Clearly, this correspondence induces a natural embedding
$\Hom \;{\bar\Mc^{\Ha^*}}_\Ha\to\Hom \;{\bar\Mc^{\Ha^*}}$.
Note
that this embedding is not an isomorphism, in general.
\begin{propn}
\label{H-H*}
The strong monoidal functor $\Hom \;{\bar\Mc^{\Ha^*}}_\Ha\to\Hom
\;{\bar\Mc^{\Ha^*}}$,
$\phi \mapsto \phi'$, given by  (\ref{eq_HH*}),
induces a strong monoidal functor ${\bar\Mc^{\Ha^*}}_\Ha\to
\bar\Mc^{\Ha^*}$.
\end{propn}
\noindent
The proof of this proposition uses the diagram technique, the properties
of
permutations $\{\tau_A\}$ and $\{\tau^A\}$, and relies on
Lemma \ref{lm_t-t}. The details are left to the reader.
\subsection{Dynamical extension of a monoidal category over a module
category}
\label{ssDE}
\label{coordinate}
  Let $\O$ be a monoidal category  and $\B$ its left module category, see
\cite{O}.
For example, $\B$ is a monoidal category and $\O$ its monoidal
subcategory.
We denote the tensor product in $\O$ and action of $\O$ on $\B$ by the
same symbol $\tp$.
For simplicity, all monoidal categories are assumed to be strict (with
trivial associativity);
the same is assumed for actions on module categories.

Let us define a \select{dynamical extension}, $\bar\O_{\tr\B}$, of $\O$
over $\B$ in the following way.
The collection of objects in $\bar \O_{\tr\B} $ coincides with that of
$\O$.
An object  $A$  of $\bar \O_{\tr\B}$ is treated as a  functor from $\B$ to
$\B$,
namely $X\stackrel{A}{\longrightarrow} A\tp X$ for all $X\in \Ob\; \B$.
Morphisms of $\bar \O_{\tr\B}$ are natural transformations of the
functors.
Namely,
$\phi\in \Hom_{\bar O_{\tr\B}}(A,B)$ is a collection $\{\phi_X\}$ of
morphisms
$\phi_X\in \Hom_\B(A\tp X,B\tp X\}$ such that
\be
\label{invariance}
\phi_X \circ(\id_A \tp \xi)=(\id_B\tp\xi)\circ\phi_{X'}
\ee
for any $\xi\in \Hom_\B(X',X)$.
The composition of morphisms in  $\bar \O_{\tr\B}$ is "pointwise",
$(\phi \bar\circ\psi)_X=\phi_X\circ\psi_X$.
Obviously, the condition (\ref{invariance}) holds for $\bar \circ$.
Clearly, $\bar O_{\tr\B}$ defined in this way is a category.
\begin{propn}
$\bar \O_{\tr\B}$ is a monoidal category with respect to the tensor product
on the objects
as in $\Ob\:\O$ and defined on the morphisms by
\be
\label{dynTP}
\begin{array}{rcll}
(\phi \bar\tp \psi)_X:=( \id_C\tp\psi_{X})\circ(\phi_{B\tp X})=(\phi_{D\tp
X})\circ (\id_A\tp\psi_{X}),
\end{array}
\ee
for $\phi\in \Hom_{\bar \O_{\tr\B}}(A,C)$, and $\psi\in \Hom_{\bar
\O_{\tr\B}}(B,D).$
\end{propn}
\begin{proof}
Let us check that the family $\{(\phi  \bar\tp \psi)_X\}$ defines  a
morphism of
functors, $A\tp B\to C\tp D$. First of all, observe that condition
(\ref{invariance}) is satisfied.
We will show that operation (\ref{dynTP}) is functorial.
Take $\{\al_X\}\in \Hom_{\bar \O_{\tr\B}}(A',A)$ and $\{\bt_X\}\in
\Hom_{\bar \O_{\tr\B}}(B',B)$.
We have for $(\phi\bar \circ\al)\bar \tp (\psi\bar \circ\bt)$:
\be
\bigl(\id_C\tp(\psi_{X}\circ\bt_X)\bigr)\circ(\phi_{B'\tp X}\circ\al_{B'\tp
X})&=&
(\id_C\tp\psi_{X})\circ(\id_C\tp \bt_X)\circ\phi_{B'\tp X}\circ\al_{B'\tp
X}
\nn\\
&=&
(\id_C\tp\psi_{X})\circ\phi_{B\tp X}\circ(\id_C\tp \bt_X)\circ\al_{B'\tp
X}
\nn\\
&=&
(\phi\bar\tp \psi)_{X}\circ(\al\bar\tp \bt)_X
\nn
\ee
for all $X\in \Ob\;\B$.
In transition to the middle line we used the condition (\ref{invariance}),
in order to permute the morphisms
$\id_C\tp \bt_X$ and $\phi_{B'\tp X}$. To prove associativity,
we take $\zeta\in\Hom_{\bar \O_{\tr\B}}(A,U)$, $\phi\in\Hom_{\bar
\O_{\tr\B}}(B,V)$,
$\psi\in\Hom_{\bar \O_{\tr\B}}(C,W)$ and find that
$\bigl(\zeta \bar\tp (\phi\bar\tp\psi)\bigr)_X$ and  $\bigl((\zeta \bar\tp
\phi)\bar\tp\psi\bigr)_X$
are equal to the same composition map
$$
A\tp B \tp C \tp X \stackrel{\zeta_{B\tp C \tp X}}{\longrightarrow}
U\tp B \tp C \tp X \stackrel{\phi_{C\tp X}}{\longrightarrow}
U\tp V \tp C \tp X\stackrel{\psi_{X}}{\longrightarrow}
U\tp V \tp W \tp X.
$$
This completes the proof.
\end{proof}
\begin{definition}
The category $\bar \O_{\tr\B}$ is called \select{dynamical extension} of
$\O$ over $\B$.
\end{definition}
\begin{remark}
Similarly to $\bar\O_{\tr \B}$, one can define a dynamical extension,
${}_{\B\tl}\!\bar\O$,
of a monoidal category $\O$ over its \select{right} module category $\B$.
Thus, the set $\Hom_{{}_{\B\tl}\!\bar\O}(A,B)$
is formed by families $\{{}_{X}\!\psi\}$ from $\Hom_{\B}(X\tp A\to X\tp
B)$
subject to the natural condition analogous to (\ref{invariance}).
The composition $\bar \circ$ is defined as the composition of functor
morphisms,
similarly to the $\bar \O_{\tr \B}$ case.
Formula (\ref{dynTP}) for
tensor products of morphisms is changed to
\be
\begin{array}{rcll}
{}_{X}\!(\phi  \bar\tp \psi)&:=&({}_{X}\!\phi\tp \id_D)\circ {}_{X\tp
A}\!\psi,& \phi\in \Hom_{{}_{\B\tl}\!\bar\O}(A,C),\quad
\psi\in \Hom_{{}_{\B\tl}\!\bar\O}(B,D).
\end{array}
\ee
\end{remark}

\subsection{Comparison of categories $\bar{\O}_{\tr \B}$ and
${}_{\B\tl}\!\bar{\O}$ with $\bar {\Mc}_{\Ha;\Lc}$
and $\bar {\Mc}^{\Ha^*}$}

Let $\Lc$ be a base algebra over a Hopf algebra $\Ha$. Let $\B$ be the
category of  left $\Lc$-modules, and $\O$ the category  $\Mc_\Ha$
of locally finite
left $\Ha$-modules.
Then $\B$ is a left $\O$-module category. The tensor product of
$A\in \Ob\;\O$ and $X\in\Ob\;\B$ is an $\Lc$-module
\be
\ell\btr (a\tp x)=  \ell^{(1)}\tr a\tp \ell^{[2]}\btr x, \quad \ell\in \Lc,
\> a\in A,\> x\in X,
\ee
where $\btr$ denotes the action of $\Lc$ and $\tr$ the action of $\Ha$.

Consider the dynamical extension $\bar\Mc_{\Ha;\Lc}$ of $\Mc_\Ha$ over the
base algebra $\Lc$ as in Example \ref{exDHL}.
Let $\psi$ be a morphism from $\Hom_{\bar \Mc_{\Ha;\Lc}}(A,B)$.
Consider the family of
maps $\psi_X\colon A\tp X \to B\tp X$, $X\in \B$, defined by the
composition
\be
A\tp X \stackrel{\psi\tp \id_X}{\longrightarrow} B\tp \Lc\tp X
\stackrel{\id_B \tp \btr}{\longrightarrow} B\tp X.
\label{eq_LX}
\ee
The maps (\ref{eq_LX}) are $\Lc$-equivariant, due to quasi-commutativity of
$\Lc$.
The following proposition is immediate.
\begin{propn}
\label{H-X}
The correspondence $\psi\mapsto \{\psi_X\}$ between morphisms induces
a strong monoidal functor $\bar \Mc_{\Ha;\Lc}\to \bar\O_{\tr \B}$ identical
on objects.
\end{propn}

Now take $\O$ to be the category $\Mc^{\Ha^*}$ of right locally finite
$\Ha^*$-comodules also considered as
left $\Ha$-modules. Put $\B$ the category of locally finite
$\Ha$-modules.
\begin{propn}
\label{OBH*}
There exists a strong monoidal functor
$\bar \Mc^{\Ha^*}\to {}_{\B\tl}\!\bar \O$.
\end{propn}
\begin{proof}
We will give a sketch of proof.
Categories ${}_{\B\tl}\!\bar \O$
and $\bar \Mc^{\Ha^*}$ have the same collection of objects, and the
functor in question is set to be
identical on objects. Let us define it on morphisms.
Let $f:\Ha^*\ot A\to B$ be a morphism in $\bar \Mc^{\Ha^*}$.
For every finite dimensional $\Ha$-module $X$ there is a natural
map $X^*\tp X\to \Ha^*$, where $X^*\tp X$ is considered as
the (left) dual module to the space of right endomorphisms (over $k$) of
$X$.
Hence, $f$ defines a collection of $\Ha$-equivariant maps $X^*\tp X\tp A\to
B$,
or, equivalently, a collection $\{f_X\}$ of
$\Ha$-equivariant maps $X\tp A \to X\tp B$.
This family extends to all locally finite
$\Ha$-modules $X$.
Thus we have built an embedding of morphisms
$\Hom\;{\bar \Mc^{\Ha^*}}\to\Hom\;{}_{\B\tl}\!\bar \O$, $f\mapsto\{f_X\}$.
It remains to check that the above correspondence is functorial
and respects the composition and the tensor product of morphisms.
We leave the details to the reader.
\end{proof}
\begin{remark}
\label{iso}
The functor from Proposition \ref{OBH*} is an isomorphism
when $\Ha^*$ decomposes into the direct sum of $X^*\tp X$, where $X$ runs
over simple $\Ha$-modules.
\end{remark}
\subsection{Dynamical associative algebras}
\label{ssDA}
The notion of dynamical associative  algebra
as an algebra in a monoidal (dynamical) category is introduced in the
standard
way.
Below we give examples of dynamical associative algebras in the categories
$\bar \Mc^{\Ha^*}$
and $\bar \Mc_{\Ha;\Lc}$.
\begin{example}[Dynamical algebras in $\bar \Mc_{\Ha;\Lc}$]
\label{exDA}
Let us consider the dynamical extension $\bar \Mc_{\Ha;\Lc}$ of the
category $\Mc_\Ha$ over
a base $\Ha$-algebra $\Lc$. An algebra $\A$ in $\bar \Mc_{\Ha;\Lc}$  is an
object
equipped with a morphism $\A\tp \A\to \A$ obeying
the associativity axiom. In terms of $\Mc_{\Ha}$, this is
equivalent to Definition \ref{defDA}. Namely,
the multiplication in $\A$ is an $\Ha$-equivariant map
$\divideontimes\colon\A\tp \A\to \A\tp \Lc$,
which is shifted associative in the sense of (\ref{shifted_assoc0}).

Now let us prove Proposition \ref{propCSP}.
This is a corollary of the following general fact.
Let $(\Cc, \tp, 1_\Cc)$ be a monoidal category whose objects are vector
spaces over $k$.
Suppose there is an object $\A\in \Ob\:\Cc$, a morphism $\iota\colon
1_\Cc\to \A$, and
an operation
$\Hom_\Cc(X,\A)\tp_k\Hom_\Cc(Y,\A)\stackrel{\circledast}{\longrightarrow}
\Hom_\Cc(X\tp Y,\A)$
for all $X,Y \in \Ob\:\Cc$.
We say that $\circledast$ is
a) natural if $(\phi\circ\al)\circledast(\psi\circ
\bt)=(\phi\circledast\psi)\circ(\al\tp \bt)$,
b) associative if
$(\phi\circledast\psi)\circledast\vartheta=\phi\circledast(\psi\circledast\vartheta)$,
and
c) unital if $\phi\circledast (\iota\circ\chi)=\phi\tp \chi$,
$(\iota\circ\chi)\circledast \phi=\chi\tp \phi$ for all morphisms $\chi$
with
 target in $1_\Cc$.
The multiplication $\mathrm{m}$ in $\A$ and the operation $\circledast$ are
related
by $\mathrm{m}=\id_\A\circledast\id_\A$, $\phi\circledast\psi =
\mathrm{m}\circ(\phi\tp \psi)$.
Now let $\A$ be an algebra in $\bar \Mc_{\Ha;\Lc}$. The unit morphism
$\iota\colon k\to \A$
in $\bar \Mc_{\Ha;\Lc}$ gives the unit map $k\to \A\tp \Lc$ in
$\Mc_{\Ha}$.
The multiplication $\divideontimes$
in $\A$ defines a natural associative unital operation on morphisms
from $\Hom\:{\bar \Mc_{\Ha;\Lc}}$ with target in $\A$.
Hence it defines a natural associative unital operation
on morphisms from $\Hom \:{\Mc_{\Ha}}$ wight target in $\A\tp \Lc$.

\end{example}

\begin{example}[Dynamical algebras in $\bar \Mc^{\Ha^*}$]
\label{algH*}
Let us describe dynamical associative  algebras in the category $\bar
\Mc^{\Ha^*}$.
The multiplication in an algebra $\A\in \Ob \; \bar \Mc^{\Ha^*}$
is an $\Ha$-equivariant map $\Lambda:{\Ha^*}\tp \A\tp \A\to \A$.
Associativity, in terms of  $\Mc^{\Ha^*}$, is formalized by the
requirement
that the following diagram is commutative:
\begin{equation}
\begin{array}{ccccccc}
\Ha^*\tp \A\tp\A\tp \A &
\stackrel{\scriptstyle\Delta }
{\longrightarrow} &\Ha^*\tp \Ha^*\tp \A \tp \A\tp \A&
\stackrel{\scriptstyle{\Lambda}}
{\longrightarrow}   & \Ha^*\tp  \A\tp\A&
\stackrel{\scriptstyle{\Lambda}}
{\longrightarrow}&  \A
\\
&&{\scriptstyle\tau^\A}
\downarrow
\hspace{35pt}
&&&&\parallel\\
&& \Ha^*\tp \A \tp \Ha^*\tp \A\tp \A&
\stackrel{\scriptstyle{\Lambda}}
{\longrightarrow}& \Ha^*\tp \A \tp \A
&
\stackrel{\scriptstyle{\Lambda}}
{\longrightarrow}&  \A
\end{array}.
\label{shifted_assoc1}
\end{equation}
This diagram is a "partial dualization" of  the diagram
(\ref{shifted_assoc0}).
The algebra $\A$ is unital if there is an element $1\in \A$ such that
$\Lambda(\la,1,a)=\Lambda(\la,a,1)=\ve(\la)a$,
for all $a\in \A$, $\la\in\Ha^*$.

The map $\Lambda$ defines a family of bilinear operations $\st{\la}$
depending on elements $\la\in \Ha^*$.
In terms of $\st{\la}$, the "shifted" associativity (\ref{shifted_assoc1})
reads (summation implicit)
\be
\label{multiplications}
(a\st{\la^{(2)}} b ) \st{\la^{(1)}} c=a^{[0]} \st{\la^{(1)}}(b\st{\la^{(2)}
a^{(1)}}c).
\ee
It is proposed by I. Kantor to consider the multiplication map  $\Ha^*\tp
\A\tp \A \stackrel{\MM}{\longrightarrow}\A $
as a ternary operation $\la\tp a\tp b\mapsto (\la a b)$ which is
associative in the sense
$$
\bigl(\la^{(1)} (\la^{(2)} a b) c\bigr)=\bigl(\la^{(1)}  a' (\la^{(2)''} b
c)\bigr),
$$
where $a'\tp \la'' = \tau^\A(\la\tp a)$, the permutation (\ref{tau*}).

\end{example}
\section{Categorical approach to quantum DYBE}
\label{sDCRA}
\subsection{Dynamical twisting cocycles}
\label{ssDC}
In this subsection we study transformations of dynamical categories.
Recall that a  functor $\tilde\Cc\stackrel{\Upsilon}{\longrightarrow} \Cc$
between two monoidal categories is called monoidal if there is
a functor isomorphism $F$ between
$\Upsilon(A) \tp \Upsilon(B)$ and $\Upsilon(A\tilde\tp B)$.
This implies a family of isomorphisms,
$$\Upsilon(A) \tp \Upsilon(B)\stackrel{F_{A,B}}{\longrightarrow}\Upsilon(A
\tilde\tp B)$$
fulfilling the cocycle conditions (for simplicity, we assume the trivial
associator)
\be
\label{cocycle}
F_{A\tilde\tp B,C}\circ \bigl(F_{A,B}\tp \id_{C}\bigr)&=&
F_{A,B\tilde\tp C}\circ \bigl(\id_{A}\tp F_{B,C}\bigr),
\\
F_{A,1}=&\id_A&=F_{1,A},
\label{norm}
\ee
where $1$ is the unit of $\Cc$.
We are mostly  interested in the situation when $\Ob\;\tilde\Cc =\Ob\;
\Cc$
and $\Upsilon$ is identical on objects.

Suppose $F$ is a cocycle in  $\Cc$, i.e. a family of invertible morphisms
$F_{A,B}\in \Aut_\Cc(A\tp B)$ fulfilling the conditions
(\ref{cocycle}) and (\ref{norm}).
Then it is possible to define a new monoidal structure on $\Cc$.
It is the same on  objects and
defined by
\be
\label{new_prod}
\phi\tilde \tp \psi :=F\circ (\phi \tp \psi)\circ F^{-1}
\ee
on morphisms.
This new monoidal category $\tilde \Cc$ coincides with the old one
if $F$ respects morphisms of  $\Cc$, i.e.
\be
\phi\tp \psi =F\circ (\phi \tp \psi)\circ F^{-1}
\ee
for all $f,g\in \Hom\;\Cc$.
\begin{remark}
\label{nonstrict}
On can define the category $\tilde \Cc$ using an arbitrary family
$F_{A,B}\in \Aut_\Cc(A,B)$ of morphisms,
which is not necessarily a cocycle. Then $\tilde \Cc$ will not be strictly
monoidal, but rather with
the associator $\Phi_{A,B,C}=F_{A,BC} F_{B,C}F_{A,B}^{-1}F_{AB,C}^{-1}$,
which satisfies the pentagon identity
in $\tilde \Cc$.
The identity functor $\tilde\Cc\to \Cc$ yields  an isomorphism of monoidal
categories.
\end{remark}
\begin{definition}[Dynamical twist]
Let $\bar \O$ be a dynamical extension of a monoidal category $\O$.
\select{Dynamical twist}
is a cocycle in $\bar \O$ that respects morphisms from $\O$.
\end{definition}
A dynamical twist is identical on $\O$, therefore $\O$ remains a
subcategory
in the twisted category $\tilde {\bar \O}$.

One of the applications of twist is a transformation of algebras.
Any cocycle $F$ in a category $\Cc$ makes a $\Cc$-algebra with the
multiplication  $\mathrm{m}$ into a $\tilde\Cc$-algebra,
with the multiplication $\mathrm{m}\circ F^{-1}$.
Let us apply this to the specific situation of dynamical twist and build
an $\bar \O$-algebra out of $\O$-algebra.
\begin{propn}
Let $F$ be a dynamical twist in $\bar \O$.
Let $\A$ be an algebra in $\O$ with multiplication
$\mathrm{m}$. Then the multiplication
$\mathrm{m}\circ F$ makes $\A$ a dynamical associative algebra, i.e. an
algebra in $\bar \O$.
\end{propn}
\begin{proof}
It follows from (\ref{new_prod}) that dynamical twist preserves $\O$ as a
monoidal subcategory
in $\tilde {\bar \O}$. Therefore $\A$ turns out to be an algebra in $\tilde
{\bar \O}$ as well.
The family $\{F^{-1}_{A,B}\}$ is a cocycle in $\tilde {\bar \O}$; the
corresponding twist of $\tilde {\bar \O}$
gives ${\bar \O}$. Applying this inverse twist to the algebra $\A$ we
obtain an ${\bar \O}$-algebra
with the multiplication $\mathrm{m}\circ F$.
\end{proof}

Below we specialize the cocycle equations (\ref{cocycle}) and (\ref{norm})
for various types of dynamical categories.

\begin{example}[Dynamical twist in $\bar O_{\tr\B}$]
Let us express a cocycle in dynamical category $\bar O_{\tr\B}$ in terms of
$\O$ and $\B$.
A cocycle in $\bar \O_{\tr\B}$ is a collection $(F_{V,W})_X$
from $\Aut_{\B}(V\tp W\tp X)$, $V,W\in \Ob\: \O$, $X\in\Ob\: \B$,
satisfying conditions
\be
(F_{V\tp W,U})_X\circ (F_{V,W})_{U\tp X} &=& (F_{V, W\tp U})_X\circ (F_{W,
U}),
\\
(F_{V,1_{\bar \O}})_X =&\id_{X\tp V}&= (F_{1_{\bar \O},V})_X.
\ee
\end{example}
\begin{example}[Drinfeld associator as a twist in ${}_{\O\tl}\bar O$]
\label{Phi-twistor}
Let  $\g$ be a complex simple Lie algebra.
In \cite{EE1}, Enriquez and Etingof proposed
a quantization of the Alekseev-Meinrenken dynamical r-matrix
using the Drinfeld associator $\Phi\in \U^{\tp 3}(\g)[[t]]$.
This quantization can be interpreted as a twist in the category
${}_{\O\tl}\bar O$, where
$\O$ be the category of free $\C[[t]]$-modules of finite rank with
$\U(\g)[[t]]$-action.
Indeed, let us put ${}_X\!(F_{A,B}):=\Phi_{X,A,B}$. Then the pentagon
identity on $\Phi$
takes the form
$$
\Phi_{A,B,C}\circ {}_X\!(F_{A\tp B,C}) \circ
{}_X\!(F_{A,B})={}_X\!(F_{A,B\tp C})\circ  {}_{X\tp A}\!(F_{B,C}).
$$
The twisted dynamical category is not strictly monoidal, cf. Remark
\ref{nonstrict}.
It is equipped with the associator $\{\Phi_{A,B,C}\}$.
\end{example}
\begin{example}[Dynamical twist in $\bar \O_\Lc$]
Consider a cocycle in $\bar \O_\Lc$, the dynamical extension of a category
$\O$ over a base algebra $(\Lc,\tau)$.
In terms of $\O$, condition (\ref{cocycle}) reads
\be
\mathrm{m}_{\Lc\tp \Lc}\circ F_{V\tp W,U}\circ (\id_{V\tp W}\tp\tau_U)\circ
F_{V,W} &=& \mathrm{m}_{\Lc\tp \Lc}\circ F_{V, W\tp U} \circ F_{W, U},
\ee
where $F_{V,W}\in \Hom_\O(V\tp W,V\tp W\tp \Lc )$ (the $\id$-isomorphisms
are dropped from the formulas).
\end{example}

\begin{example}[Dynamical twist in $\bar \Mc^{\Ha^*}$]
\label{exCOH*}
Let us specialize the notion of cocycle for
the category $\bar \Mc^{\Ha^*}$.
A morphism $\Ha^*\tp A\stackrel{f}{\to} B$ in the  category $\Mc^{\Ha^*}$
can be thought of as a family
of maps $f^\la\colon A\to B$ parameterized by elements $\la\in \Ha^*$.
Let $\Omega^\la$ be a family of linear operators on the tensor product
$\tp_{l=1}^m V_l$
of $\Ha^*$-comodules $V_l$, $l=1,\ldots, m$.
By $^{V_i}\!\Omega^\la$, or simply by $^i\!\Omega^\la$,
we denote a family of linear operators on $\tp_{l=1}^m V_l$ defined by
$$
{\: }^i\! \Omega^\la(v_1\tp \ldots \tp v_m) := \Omega^{\la
v_i^{(1)}}(v_1\tp \ldots \tp v_i^{[0]}\tp \ldots \tp v_m),
$$
where $v_i^{[0]}\tp v_i^{(1)}$ denotes the right $\Ha^*$-coaction
$\delta(v_i)$ (as always, the summation is implicit).
The collection of morphisms $F^{\la}_{V,W}\in \Hom_{\Mc^{\Ha^*}}({\Ha^*}\tp
V\tp W, V\tp W)$ satisfies
condition (\ref{cocycle}) and (\ref{norm}) in $\bar\Mc^{\Ha^*}$ if and only
if
\be
\label{twistH*}
F^{\la^{(1)}}_{V\tp W, U}F^{\la^{(2)}}_{V,W}
&=&
 F^{\la^{(1)}}_{V, W \tp U}{\:}^V\!F^{\la^{(2)}}_{W,U},
\\
F^\la_{V,k} =&\id_V&= F^\la_{k,V}.
\ee
\end{example}
\begin{example}[Universal cocycle]
\label{exUC}
Assume that $\Ha$
is a Hopf subalgebra of another Hopf algebra, $\U$.
Then the category $\Mc_{\U}$ is a subcategory of the category
$\Mc_{\Ha}$. Let $\Lc$ be a base algebra over $\Ha$.
Suppose there is an invertible element
$\bar\F=\bar\F_1\tp\bar\F_2\tp\bar\F_3\in \U\tp \U\tp\Lc $
that satisfies the condition
\be\label{coninv1}
h^{(1)}\bar\F_1\tp h^{(2)}\bar\F_2\tp h^{(3)}\tr
\bar\F_3=\bar\F_1 h^{(1)}\tp \bar\F_2 h^{(2)}\tp \bar\F_3
\ee
for all $h\in \Ha$, and the conditions
\be
\label{shifted_cocycle}
(\Delta\tp \id)(\bar\F)\;{}^{(3)}\!\bar\F_{12}
&=&
(\id\tp \Delta)(\bar\F)(\bar\F_{23}),
\\
(\ve \tp \id\tp \id)(\bar\F)=&1\tp 1\tp 1& = (\id\tp \ve \tp \id)(\bar\F)
\ee
in $\U\tp\U\tp \U\tp \Lc$. Here notation ${}^{(3)}\!\bar\F$ means
$\delta(\bar\F)$, where $\delta$ is the coaction $\Lc\to \Ha\tp \Lc$; the
$\Ha$-component
is embedded to the third tensor factor in  $\U\tp\U\tp \U\tp \Lc$.  The
element $\bar\F$ defines a cocycle
in $\bar \Mc_{\U;\Lc}$:
for $\U$-modules $V$ and $W$ one has $F_{V,W}:=\rho_V(\bar\F_1)\tp
\rho_W(\bar\F_2)\tp\bar\F_3$.
This cocycle clearly respects morphisms in $\Mc_\U$, hence it is a
dynamical twist.
The element $\bar\F$ may be called a \select{universal dynamical twist}, by
the analogy with
universal R-matrix.
Equation (\ref{shifted_cocycle}) turns to the shifted cocycle condition of
\cite{Xu2} for
$\Ha$ being a universal enveloping algebra.
\end{example}
\subsection{Quantum dynamical R-matrix}
\label{ssDB}
\subsubsection{Dynamical Yang-Baxter equation}
Let us consider the notion of the Yang-Baxter equation in
dynamical categories. Let $\Cc$ be a braided monoidal category
with braiding $\si$. The braiding is a collection, $\{\si_{A,B}\}$, of
morphisms
$A\tp B\stackrel{\si_{A,B}}{\longrightarrow} B\tp A$
for $A,B\in \Ob \; \Cc$
obeying conditions
\be
\si_{A,B}\circ\si_{A,C}\circ\si_{B,C}
&=&
\si_{B,C}\circ\si_{A,C}\circ \si_{A,B},
\label{YBE}
\\
\si_{A\tp B,C}=\si_{A,C}\circ \si_{B,C},&&\si_{C,A\tp B}=\si_{C,B}\circ
\si_{C,A}
\label{hex}
\ee
and respecting morphisms, i.e. $(f\tp g)\circ\si =\si \circ(g\tp f)$
for all $f,g\in\Hom\: \Cc$
(in fact, (\ref{YBE}) follows from  (\ref{hex}) and functoriality of
$\si$).
Condition (\ref{YBE}) is called the
Yang-Baxter equation, conditions  (\ref{hex}) are, in fact, the hexagon
identities. If $\si$ fulfils (\ref{YBE}) and (\ref{hex})
but is not functorial (does not respect morphisms), we call it
\select{pre-braiding}. This is the case
when $\si$ is a braiding in a subcategory $\Cc'$ of $\Cc$
such that $\Ob \;\Cc'=\Ob \;\Cc$, e.g.,
when $\Cc$ is a dynamical extension of a $\Cc'$.
Then $\Cc$ has more morphisms than $\Cc'$, and they are not respected by
$\si$,
in general\footnote{For instance, the dynamical extension $\bar\O_\Lc$
of a braided category $(\O,\si)$ over a commutative algebra $\Lc$ in $\O$
(cf. Example \ref{qcomalg}) is
braided if and only if $\si_{A,\Lc}\circ\si_{\Lc,A}=\id_{\Lc\tp A}$
for all $A\in \Ob \; \O$, i.e.,  when  $\O$ is a symmetric
category.}.

Given a pre-braiding $\si$ in $\Cc$, it is possible to restrict it to
a braiding in a subcategory $\Cc_\si$ defined as follows.
The objects in $\Cc_\si$ are those of $\Cc$.
A morphism $f\in\Hom_{\Cc}(A,B)$ is a morphism from $\Hom_{\Cc_\si}(A,B)$
if and only if
\be
\si_{B,C}\circ (f\tp \id_C)=(\id_C\tp f)\circ\si_{A,C}
,\quad
(f\tp \id_C)\circ\si_{C,A}=\si_{C,B}\circ(\id_C\tp f)
\ee
for all $C\in \Ob\;\Cc$.

\begin{propn}
$\Cc_\si$ is a braided category with braiding $\si$.
\end{propn}
\begin{proof}
It follows from (\ref{YBE}) and (\ref{hex}) that $\si$ lies in $\Cc_\si$.
Condition (\ref{hex}) guarantees that $\Cc_\si$ is a monoidal
category. Therefore, $\si$ is a pre-braiding in $\Cc_\si$ and respects
morphisms
in it by construction; hence $\si$ is a braiding in $\Cc_\si$.
\end{proof}
\begin{propn}
\label{twYBE}
Let $\si$ be a pre-braiding in $\Cc$ and let $F$ be a cocycle in $\Cc$
respecting
morphisms from $\Cc_\si$.
Then the family
\be
\label{tw_brad}
\bar\si_{A,B}:=F^{-1}_{B,A}\circ\si_{A,B}\circ F_{A,B}
\ee
satisfies the Yang-Baxter equation (\ref{YBE}).
\end{propn}
\begin{proof}
Define $\Omega_{A,B,C}:=F_{A\tp B,C}\circ(F_{A,B}\tp \id_C)=F_{A,B\tp
C}\circ(\id_A \tp F_{B,C})$
for all $A,B,C\in \Ob\: \Cc$.
Since $F$ respects morphisms from $\Cc_\si$,
we have $\Omega^{-1}_{A,C,B}\circ(\id_A\tp \si_{B,C})\circ
\Omega_{A,B,C}=\id_A\tp \bar \si_{B,C}$
and $\Omega^{-1}_{B,A,C}\circ(\si_{A,B}\tp \id_C)\circ \Omega_{A,B,C}=\bar
\si_{A,B}\tp \id_C$
for all $A,B,C$.
Multiplying the equation (\ref{YBE}) by $\Omega_{C,B,A}^{-1}$ from the
left
and by $\Omega_{A,B,C}$ from the right, we prove the statement.
\end{proof}
Applied to dynamical twists, the Proposition \ref{twYBE} yields the
following corollary.
\begin{corollary}
Let $\O$ be a braided category with the braiding $\si$. Let $\bar \O$ be
a dynamical extension of $\O$ and $F$ a dynamical twist in $\bar \O$.
The collection of morphisms (\ref{tw_brad}) for $A,B\in \bar \O$,
satisfies the Yang-Baxter equation in $\bar \O$.
\end{corollary}

In general, a twist destroys the hexagon identities in the twisted category
$\tilde \Cc$.
However, it yields a pre-braiding in an equivalent category to  $\tilde
\Cc$, which is constructed
in Subsection \ref{sssDB}.
Another way to fix the situation is  when $\Cc=\bar \Mc_{\Ha,\Lc}$, the
dynamical extension of the category of $\Ha$ modules
over a base algebra $\Lc$. Namely, there is a realization of  $\bar
\Mc_{\Ha,\Lc}$
as a category of modules  over a certain bialgebroid, \cite{DM5}.  A
dynamical twist gives rise to a bialgebroid twist,
which transforms the braiding in the category of modules over the
bialgebroid.

We call a solution of (\ref{YBE}) a \select{dynamical R-matrix}.

Below we specialize this definition of dynamical R-matrix to various types
of dynamical categories.

\begin{example}[Dynamical R-matrix in $\bar \O_{\tr\B}$]
The dynamical R-matrix in the category $\bar \O_{\tr\B}$ is defined by
\be
(\si_{A,B})_X\circ(\si_{A,C})_{B\tp
X}\circ(\si_{B,C})_X&=&(\si_{B,C})_{A\tp
X}\circ(\si_{A,C})_X\circ(\si_{A,B})_{C\tp X},
\ee
where $\si$ is a collection of invertible morphisms $(\si_{A,B})_X\in
\Aut_{\B}(A\tp B\tp X)$.
\end{example}
\begin{example}[Dynamical R-matrix in $\bar \O_\Lc$]
Consider the category $\bar \O_\Lc$, a dynamical extension of a monoidal
category $\O$ over a base algebra $(\Lc,\tau)$,
cf. Subsection \ref{ssExtL}. Let $\mathrm{m}$ be the multiplication in the
algebra $\Lc$ and $\mathrm{m}^3$ denote the three-fold product
$\mathrm{m}\circ (\mathrm{m}\tp \id_\Lc)$.
In terms of $\O$ and $(\Lc,\tau)$, equation (\ref{YBE}) reads
\be
\mathrm{m}^3\circ\si_{A,B}\circ\tau_B\circ\si_{A,C}\circ\si_{B,C}
&=&
\mathrm{m}^3\circ\tau_A\circ\si_{B,C}\circ\si_{A,C}\circ\tau_C\circ\si_{A,B}
\ee
where $\si_{A,B}\in \Hom_{\hat \O}(A\tp B,A\tp B\tp \Lc )$.
\end{example}
\begin{example}[Dynamical R-matrix in $\bar\Mc^{\Ha^*}$]
We use the notation of Example \ref{exCOH*}.
A collection of morphisms $\{\si^{\la}_{A,B}\}$ from $\Hom\:{\Mc^{\Ha^*}}$
fulfills the equation
(\ref{YBE}) in $\bar\Mc^{\Ha^*}$ if and only if
\be
^A\!\si^{\la^{(1)}}_{B,C}\>\>\si^{\la^{(2)}}_{A,C}\>\> {\:
}^C\!\si^{\la^{(3)}}_{A,B}
&=&
\si^{\la^{(1)}}_{A,B} \>\>{\: }^B\!\si^{\la^{(2)}}_{A,C}
\>\>\si^{\la^{(3)}}_{B,C}.
\ee
\end{example}
\begin{example}[Universal dynamical R-matrix]
\label{exUR} Consider the situation of Example \ref{exUC} assuming that
$\Ha$ is a Hopf subalgebra in another Hopf algebra, $\U$, and
$\Lc$ is a base $\Ha$-algebra. Definition \ref{def_QDR} of
Subsection \ref{ssDYBEs} introduces a universal quantum dynamical
R-matrix of $\U$ over the base $\Lc$.  For any pair $V$ and $W$ of
$\U$-modules considered as modules over $\Ha$, it gives
$\si_{V,W}:=P_{V,W}R_{V,W}$, where $P$ is the usual flip and
$R_{V,W}=(\rho_V\tp \rho_W)(\bar R)$ is the image of $\bar R$ in
$\End(V)\tp \End(W)\tp \Lc$.
\begin{propn}
Suppose the Hopf algebra $\U$ is quasitriangular and let $\Ru$ be
its universal R-matrix. Let  $\Lc$ be a base algebra over
$\Ha\subset \U$ and $\bar\F\in \U\tp \U\tp  \Lc$ a universal
dynamical twist. Then the element $\bar  \Ru:=\bar\F_{21}^{-1}\Ru
\bar\F$ is a universal dynamical R-matrix.
\end{propn}
\begin{proof}
This statement can be checked directly. Also, it can be verified by
passing to representations of $\U$. Then it follows from
Proposition (\ref{tw_brad}).
\end{proof}
\end{example}

\subsubsection{Dynamical (pre-) braiding}
\label{sssDB}
Let $\Cc$ be a monoidal category. Let $F$ be a cocycle in $\Cc$ and
$\tilde \Cc$ be the twisted category defined in Subsection \ref{ssDC}.
Suppose $\si$ is a pre-braiding in $\Cc$. As was mentioned above, the
hexagon identities  (\ref{hex}) are destroyed in $\tilde \Cc$. We are
going to construct
an equivalent monoidal category $F(\Cc)$ where the twist of $\si$
will be a pre-braiding.

We consider formal  sequences (or words) $\Ab:=(A_1,A_2,..., A_n)$, $n>0$,
of objects
from $\Cc$. Given two words $\Ab$ and $\Bb$, let $\Ab\bullet \Bb$ denote
the
concatenation $(A_1,A_2,..., A_n,B_1,B_2,..., B_m)$.

Let $\al(\Ab)$ denote the tensor product $A_1\tp \ldots \tp A_n\in
\Ob\;\Cc$.
By induction on the length of words, let us
introduce an isomorphism $\Omega_{\Ab}$ of $\al(\Ab)\in \Ob\;\Cc$ setting
$\Omega_{\Ab}:=\id_{A}$  for $\Ab=A\in \Ob\:\Cc$ and
\be
\Omega_{\Ab\bullet \Bb}:=
   F_{\al(\Ab),\al(\Bb)}(\Omega_{\Ab}\tp \Omega_{\Bb}).
\ee
One can check, using the cocycle condition (\ref{cocycle}), that
$\Omega_{\Ab}$
does not depend on representation of $\A$ as the concatenation.
Using the family $\{\Omega_\Ab\}$, define a transformation $F(f)$ of
morphisms
$\al(\Ab) \stackrel{f}{\to} \al(\Bb)$ in $\Cc$
setting
\be
\label{F(f)}
F(f):={\Omega_{\Bb}} f {\Omega_{\Ab}}^{-1}
\ee

Let us construct the category $F(\Cc)$. The objects of $F(\Cc)$ are
finite formal sequences of objects from $\Cc$. The space of morphisms
$\Hom_{F(\Cc)}(\Ab , \Bb)$ consists of $F(f)$, where $f$
is a morphism from $\Hom_\Cc\bigl(\al(\Ab), \al(\Bb)\bigr)$.

We define the tensor product of objects $\Ab$ and $\Bb$ of $F(\Cc)$
as the concatenation  $\mathrm{\bf A}\bullet \mathrm{\bf B}$.
The empty word plays the role of the unit object.

Let us define the tensor product of morphisms in $F(\Cc)$. Let
$F(f)\colon  \mathrm{\bf A}\to \mathrm{\bf A}'$ and $F(g)\colon \mathrm{\bf
B} \to \mathrm{\bf B}'$
be two morphisms. Then we put
\be
\label{fucnt_hom}
F(f)\bullet F(g):=F(f\ot g)\colon \mathrm{\bf A}\bullet  \mathrm{\bf B}\to
\mathrm{\bf A}'\bullet  \mathrm{\bf B}'.
\ee
The category $F(\Cc)$ is equivalent to $\Cc$. Indeed,
the correspondence
$\mathrm{\bf A}\mapsto \al(\Ab)$, $F(f)\mapsto f$
gives a strong monoidal functor $\alpha\colon F(\Cc)\to \Cc$.
Consider also the functor
$\beta\colon \Cc\to F(\Cc)$ defined  on objects
by $\beta(A)=A$, the word of length $n=1$, and  on morphisms by
$\beta(f)=f$.
This functor is monoidal. Indeed, one can interpret
$\Omega_{\Ab}$ as a morphism in $F(\Cc)$, namely,
$$
\Omega^{-1}_{\Ab}\in \Hom_{F(\Cc)}(A_1\bullet ...\bullet  A_n,{A_1\tp
...\tp  A_n}  ).$$
So we obtain the transformation of the tensor products
$\bt(A)\bullet \bt(B)\stackrel{\Omega^{-1}_{(A,B)}}{\longrightarrow} \bt
(A\tp B)$.
The  functors $\al$ and $\bt$
give the equivalence of categories $\Cc$ and $F(\Cc)$.

\begin{propn}
Let $\si$ be a pre-braiding in $\Cc$.
Then the collection $\si_{\mathrm{\bf A},\mathrm{\bf
B}}:=F(\si_{\al(\mathrm{\bf A}),\al(\mathrm{\bf B})})$
is a pre-braiding in $F(\Cc)$.
\end{propn}
\begin{proof}
Apply the functor $F$ to equations
(\ref{YBE}) and (\ref{hex}) and use  the definition (\ref{fucnt_hom}).
\end{proof}
For example, let us specialize $\si_{\Ab,\Bb}$ for $\Ab=A$ and $\Bb=B$. In
this case, we have
$\Omega_{\Ab,\Bb}=F_{A,B}$. Applying formula (\ref{fucnt_hom}),
we obtain
$\si_{\Ab,\Bb}=F_{B,A}\si_{A,B}{F^{-1}}_{A,B}$ for $\Ab=A$ and $\Bb=B$.
\section{A construction of dynamical twisting cocycles}
\label{sDCC}
\subsection{Associative operations on morphisms and twists}

Let $\Cc$ be a monoidal category and $\Cc'$ a subcategory in $\Cc$.
We are going to show that cocycles in $\Cc'$ (see Subsection \ref{ssDC})
are in one-to-one correspondence with natural associative operations on
morphisms
$\Hom_{\Cc}(A,V)$, where $A \in \Ob\:\Cc$ and $V\in \Ob\: \Cc'$.
First of all observe that a cocycle $F$ in $\Cc'$ defines
such an operation by the formula $\phi\circledast
\psi:=F\circ(\phi\tp\psi)$.
The converse is also true.
\begin{lemma}
\label{tilde_tp}
Suppose there is an associative operation
$$
\Hom_{\Cc}(A,V)\tp \Hom_{\Cc}(B,W)\stackrel{\circledast}{\longrightarrow}
\Hom_{\Cc}(A\tp B,V\tp W)
$$
for all $A,B\in \Ob\: \Cc$ and $V,W\in \Ob\: \Cc'$ that is natural with
respect to its $\Cc$-arguments:
\be
\label{nat}
(\phi\circ \al)\circledast(\psi\circ \bt)=
(\phi \circledast \psi)\circ (\al \tp \bt),
\ee
whenever $\phi \in \Hom_{\Cc}(A,V)$, $\psi \in \Hom_{\Cc}(B,W)$,
$\al, \bt  \in \Hom\: {\Cc}$.
Suppose it is unital, i.e.
$$
\phi\circledast \chi = \phi\tp \chi, \quad \chi\circledast \phi = \chi \tp
\phi
$$
for any morphism $\phi$ and any $\chi \in \Hom_{\Cc}(B,1_\Cc)$.
Then the family
\be
\label{F}
F_{V,W}:=\id_{V}\circledast \id_{W}\in \End_{\Cc}(V\tp W)
\ee
is a cocycle in $\Cc'$. This cocycle respects morphisms
from a subcategory $\Cc''$ in $\Cc'$ if and only if
the operation $\circledast$ is natural with respect to $\Cc''$-arguments,
i.e.
$$
(\zeta\circ\phi)\circledast(\eta\circ\psi)=
(\zeta \tp \eta)\circ(\phi \circledast \psi)
$$
whenever $\phi \in \Hom_{\Cc}(A,V)$, $\psi \in \Hom_{\Cc}(B,W)$,
$\zeta, \eta \in \Hom\:{\Cc''}$.
\end{lemma}
\begin{proof}
By the definition (\ref{F}), the expression $F_{U\tp V,W}\circ (F_{U,V}\tp
\id_{W})$
is equal to
\be
\label{3id}
\bigl(\id_{U\tp V}\circledast \id_{W}\bigr)\circ\bigl((\id_{U}\circledast
\id_{V})\tp \id_{W}\bigr)
=\id_{U}\circledast \id_{V} \circledast \id_{W}.
\ee
Here we have used condition (\ref{nat}).
Similarly, the expression  $F_{U, V\tp W}\circ (\id_{U}\tp F_{V,W})$ is
brought
to the right-hand side of (\ref{3id}). This proves the cocycle
condition.
\end{proof}
\noindent
\subsection {Dynamical adjoint functors}
\label{ssDCF}
In this subsection we formulate the notion of dynamical adjoint functor,
which appears
to be very useful in constructing dynamical twists.
Let $\O$ be a monoidal category and $\O'$ its monoidal subcategory; the
embedding
functor is denoted by $\R$. Let $\B$ and $\B'$ be right module categories
over $\O$ and
$\O'$, respectively.
\begin{definition}
\label{def_DC}
A functor  $\B\stackrel{\M}{\longrightarrow}\B'$
is called \select{dynamical adjoint} to $\R$ if
there is an isomorphism of the following three-functors from
$\B\times\B\times\O'$ to the category of
linear spaces:
\be
 Y\!\times\! X\!\times\! V\! \to \Hom_\B\bigl(Y, X\tp \R(V)\bigr) &\simeq&
 Y\!\times\! X\!\times\! V\! \to \Hom_{\B'}\bigl(\M(Y), \M(X)\tp V\bigr).
\label{d-adjoint}
\ee
\end{definition}

Given a pair of dynamical adjoint functors, we
define an operation $\circledast$ on morphisms
from $ \Ob\:{}_{\B\tl}\!\bar\O$ to its subcategory
$\Ob\:{}_{\B\tl}\!\bar\O'$ in the following way.
A pair $\{{}_{X}\!\phi\}\in \Hom_{{}_{\B\tl}\!\bar \O}\bigl( A,\R(V)\bigr)$
and $\{{}_{X}\!\psi\}\in \Hom_{{}_{\B\tl}\!\bar \O}\bigl(B,\R(W)\bigr)$
defines a family of ${\B'}$-morphisms $\M(X\tp A\tp B)\to \M(X)\tp V\tp
W$,
for all $X\in \B$, via the composition
\be
\label{M-M}
\M(X\tp A\tp B)\stackrel{{}_{(X\tp A)}\!\tilde\psi}{\longrightarrow}
\M(X\tp A)\tp W
\stackrel{{}_{X}\!\tilde\phi\tp \id_W}{\longrightarrow} \M(X)\tp V\tp W.
\ee
By the tilde we denote the image of a morphism from $\Hom \;
{}_{\B\tl}\!\bar \O$ under the correspondence (\ref{d-adjoint}).
By condition (\ref{d-adjoint}), the  composition (\ref{M-M}) defines a
morphism,
\be
\label{assop}
X\tp A\tp B\stackrel{{}_{X}\!(\phi\circledast \psi)}{\longrightarrow} X\tp
\R(V\tp W),
\ee
in the category $\B$. Functoriality with respect to the first  argument in
(\ref{d-adjoint}) implies that the family $\{{}_{X}\!(\phi\circledast
\psi)\}$
is in fact an ${}_{\B\tl}\!\bar \O$-morphism. The associativity of the
operation
$\circledast$ follows from the associativity of composition of morphisms
in the category $\B'$.
\begin{propn}
\label{thDCF_tw_pr}
A pair of dynamical adjoint functors defines, by formula (\ref{assop}),
an associative operation $\phi\tp\psi\to \phi\circledast\psi$ that
satisfies  the
conditions of Lemma \ref{tilde_tp} for $\Cc={}_{\B\tl}\!\bar \O$ and
$\Cc'={}_{\B\tl}\!\bar \O'$. It is  $\O'$-functorial and thus yields
a dynamical twist of $\O'$.
\end{propn}
In the next subsection, using Lemma \ref{tilde_tp} and Proposition
\ref{thDCF_tw_pr},
we construct a dynamical cocycle
in the category of $\g$-modules considered as a subcategory of
$\l$-modules,
where $\l$ is an arbitrary Levi subalgebra in $\g$.

\subsection{Generalized Verma modules}
Let $\g$ be a complex reductive Lie algebra with the Cartan
subalgebra $\h$ and $\g=\n^-\oplus\h\oplus\n^+$ its polarization
with respect to $\h$.

We fix a Levi subalgebra $\l$, which is, by definition, the
centralizer of an element in $\h$. The algebra $\l$ is reductive,
so it is decomposed into the direct sum of its center and the
semisimple part, $\l=\c\oplus \l_0$, where $\l_0=[\l,\l]$.
Also, there exists a decomposition
\be
 \label{levi}
\g=\n_\l^-\oplus\l\oplus\n_\l^+,
\ee
where $\n_\l^\pm$ are subalgebras in $\n^\pm$.
Let $\p^\pm$ denote the parabolic subalgebras $\l\oplus\n_\l^\pm$.

Let $X$ be a finite dimensional semisimple representation of $\l$.
We consider $X$ as a left $\U(\l)$-module.
Being extended  by the trivial action of $\n_\l^+$ on $X$,
this  representation can be considered as a left $\U(\p^+)$-module.
We denote by $M_X$ the generalized Verma module,
$M_X:=\U(\g)\ot_{\U(\p^+)}X$. It is a left $\U(\g)$-module, and
the natural map $\U(\n_\l^-)\ot_\C X\to \U(\g)\ot_{\U(\p^+)}X$ is
an isomorphism of vector spaces.

Let us consider the dual representation $X^*$ as a left $\U(\l)$-module
with the action
\be
\label{dualact}
(u\ff)(x)=\ff(\gamma(u)x),
\ee
where $\ff\in X^*$, $x\in X$, $u\in\l$,
and $\gamma$ denotes the antipode in $\U(\g)$.
Analogously to $M_X$, we define the generalized Verma module
$M_{X^*}^-:=\U(\g)\ot_{\U(\p^-)}X^*$ naturally isomorphic as a vector
space
to $\U(\n_\l^+)\ot_\C X^*$.

There exists the following equivariant pairing between $M_{X^*}^-$
and $M_X$. Let $u_1\ot\ff\in\U(\n_\l^+)\ot_\C X^*$,
$u_2\ot x\in \U(\n_\l^-)\ot_\C X$. We put
$\langle u_1\ot\ff,u_2\ot x\rangle=\ff\bigl(s(\gm(u_1)u_2)x\bigr)$,
where $s$ is the projection $\U(\g)\to\U(\l)$ along the
direct sum decomposition
$$\U(\g)=\U(\l)\oplus(\n_\l^-\U(\g)+\U(\g)\n_\l^+).$$
It is obvious that this pairing
defines the $\U(\g)$-equivariant map
\be
\label{pairingmap}
M_{X^*}^-\to M_X^*,
\ee
where $M_X^*$ denotes
the restricted dual $\U(\g)$-module to $M_X$, which is defined
as follows. It is clear that $M_X=\oplus_\mu M_X[\mu]$,
where $M_X[\mu]$ is the finite dimensional subspace of weight
$\mu\in\h^*$. We put $M_X^*:=\oplus_\mu (M_X[\mu])^*$ with the
$\U(\g)$-action
similar to (\ref{dualact}).
It is known that
map (\ref{pairingmap}) is an isomorphism for representations $X$
satisfying Proposition \ref{Shap} below.

Since $\U(\l)=\U(\l_0)\ot\U(\c)$,
where $\l_0$ is the semisimple part of $\l$ and $\c$ its center,
a representation $X$ is an irreducible $\U(\l)$-module if and only if it
can be presented as the tensor product of two representations:
\be
\label{tepr}
X=X_0\ot\C_\lambda.
\ee
Here $X_0$ is an irreducible representation of $\l_0$, and $\C_\lambda$
is a one dimensional representation of $\c$ defined by a character
$\lambda\in \c^*$; both $X_0$ and $\C_\lambda$ are lifted to
$\U(\l)$-modules in the natural way. It is clear that representation
(\ref{tepr}) is unique. We call the element $\lambda$ from (\ref{tepr})
the \select{character} of $X$.

Let $\alpha_i$, $i=1,...,\dim \c$, be the simple roots
with respect to $\h$ that are not roots of $\l$, and
$e_{\pm\alpha_i}$ the corresponding root vectors
such that $(e_{\alpha_i},e_{-\alpha_i})=1$
for the Killing form $(.,.)$ in $\g$. Then
$h_i:=[e_{\alpha_i},e_{-\alpha_i}]$, $i=1,...,\dim \c$,
form a basis in $\c$. Denote by ${\mathcal Y}$ the union of hyperplanes
in $\c^*$ consisting of $\lambda\in\c^*$ having at least one coordinate
$\la(h_i)$ integer.
\begin{propn}[\cite{J}]
\label{Shap}
Let $X$ be a semisimple representation of $\l$. If the characters
of its irreducible components do not belong to ${\mathcal Y}$, then the
map
(\ref{pairingmap}) is an isomorphism.
\end{propn}
We call an $\l$-module $X$ generic if it satisfies this proposition.
\subsection{Dynamical twist via generalized Verma modules}
\label{ssGVM}
In this subsection we construct a dynamical cocycle for the case when
the Hopf algebra $\Ha$ is a (quantum) universal enveloping
algebra of a Levi subalgebra $\l$ in a reductive Lie algebra, $\g$.
Our method is a generalization to noncommutative and non-cocommutative Hopf
algebras
of the construction of Etingof and Varchenko, \cite{EV3}.
For simplicity we consider only classical universal enveloping algebras
$\U=\U(\g)$,  $\Ha=\U(\l)$. The construction carries over  to the quantum
groups in
the straightforward way.
Recall that $\Mc_{\U(\l)}$ and $\Mc_{\U(\g)}$ denote
the categories of locally finite semisimple
modules over $\U(\l)$ and $\U(\g)$, respectively.
\begin{lemma}
\label{M-L}
For all $Y\in \Mc_{\U(\l)}$, $V\in \Mc_{\U(\g)}$,  and generic $X\in
\Mc_{\U(\l)}$
\be
\label{Verma}
\Hom_\g(M_Y,M_X\tp V)\simeq \Hom_\l(Y,X\tp V).
\ee
\end{lemma}
\begin{proof}
Since, by Proposition \ref{Shap}, the module $M_X^*$ is
isomorphic to $M_{X^*}^-$ for generic $X$,
we have
\be
\label{eqqq1}
\Hom_\g(M_Y,M_X\tp V)\simeq \Hom_{\g}(M_Y\ot M_X^*,V)\simeq
\Hom_{\g}(M_Y\ot M^-_{X^*}, V),
\ee
where  $M^*_X$ is the restricted dual to $M_X$.
Since $M_Y\ot M^-_{X^*}\simeq\mathrm{Ind}^\g_{\l}(Y\ot X^*)$ as a
$\g$-module,
we can apply
the Frobenius reciprocity and obtain
\be
\label{eqqq2}
\Hom_{\g}(M_Y\ot M^-_{X^*}, V)\simeq \Hom_{\l}(Y\ot X^*, V )\simeq
\Hom_{\l}(Y, X\tp V).
\ee
Combining (\ref{eqqq1}) and (\ref{eqqq2}) we obtain the lemma.
\end{proof}
Set, in terms of Definition \ref{def_DC}, $\O$ to be the full subcategory
in $\Mc_{\U(\l)}$
of modules whose characters belong to the weight lattice of $\g$ relative
to $\h$.
This category contains $\Mc_{\U(\g)}$
as a subcategory, which we put to be $\O'$.
Let $\B$ be the full subcategory in $\Mc_{\U(\l)}$ of modules whose
characters do not belong to ${\mathcal Y}$;
it is a module category over $\O$.
Let $\B'$ be the category of all $\U(\g)$-modules.
Put $\R\colon \O'\to\O$ to be the restriction functor making an
$\U(\g)$-module a module over  $\U(\l)$.
We define the adjoint functor $\M$ as follows. For $X\in \Ob\Mc_{\U(\l)}$
we put $\M(X)=M_X$, the generalized Verma module corresponding to $X$.
It is clear that any morphism $X\to Y$ of $\U(\l)$-modules naturally
corresponds
to a morphism $M_X\to M_Y$ in the category $\B'$.

\begin{corollary}
\label{crlVerma_DCF}
The functor $X\stackrel{\M}{\to} M_X$ is dynamical adjoint
to the restriction functor $\Mc_{\U(\g)}\stackrel{\R}{\to}  \O$.
\end{corollary}
\begin{proof}
All we have to check is that correspondence (\ref{Verma})
is natural with respect to $Y$, $V$, and generic $X$.
This holds because the  Frobenius reciprocity gives a natural
isomorphism between adjoint functors
for generic $X$.
\end{proof}

Let us consider the category $\Mc^{\U^*(\g)}$  of locally finite semisimple
right $\U^*(\g)$-comodules.
Note that $\Mc^{\U^*(\g)}$ is naturally isomorphic to the category
$\Mc_{\U(\g)}$ of locally finite left $\U(\g)$-modules
and hence to a subcategory of locally finite semisimple left
$\U(\l)$-modules.
We call the \select{dynamical extension of} $\Mc^{\U^*(\g)}$
\select{within} $\bar\Mc^{\U^*(\l)}$
the full subcategory in $\bar\Mc^{\U^*(\l)}$ whose objects belong to
$\Mc^{\U^*(\g)}$.

Let us consider in more detail the structure of $\U^*(\l)$. First of all,
$\U^*(\l)$ can be
interpreted as the algebra of polynomial functions on the connected simply
connected
Lie group $\hat H$ corresponding to the Lie algebra $\l$. That is,
$\U^*(\l)$ is generated
over $\C$ by matrix elements of all finite dimensional representations of
$\l$.
The group $\hat H$ is presented as the Cartesian product $\hat H_0\times
\c$
of the semisimple subgroup $\hat H_0$ and $\c$ viewed as an abelian group.

It is well known that $\U^*(\l)=\oplus_V \End^*_\C(V)$, where $V$ runs over
the
irreducible $\U(\l)$-modules.
Each irreducible representation of $\l$ has the form $V=V_0\tp \C_\mu$,
where $V_0$ is a module over the semisimple part $\l_0$ of $\l$ and
$\C_\mu$ is
a one dimensional representation of the center $\c$.
Let $e^\mu\colon\c\to  \C^\times $ be the matrix element of $\C_\mu$ and
$e^{\V_0}_{ij}$ the matrix elements
of $V_0$. The elements $\{e^{V_0}_{ij} e^\mu\}$ form a basis in the vector
space $\U^*(\l)$.
We call an element $\la\in \U^*(\l)$ generic if the decomposition of $\la$
by this basis contains no $e^\mu$, where $\mu\in \mathcal{Y}$.

\begin{thm}
\label{H*Utw}
Let $\l$ be a Levi subalgebra in a reductive Lie algebra $\g$, $\U(\l)$ the
corresponding
Hopf subalgebra in the universal enveloping algebra $\U(\g)$.
There exists an $\U(\l)$-equivariant map
$$\bar\F\colon\U^*(\l)\to \U(\g)\tp \U(\g)$$
such that for generic $\la\in \U^*(\l)$ the family
$F^\la_{V,W}:= (\rho_V\tp \rho_W)\bigl(\bar\F(\la)\bigr)$,
$ V, W \in \Ob \;\Mc^{\U^*(\g)}$,
is a dynamical twist (\ref{twistH*}) in the dynamical extension of the
category $\Mc^{\U^*(\g)}$
within $\bar \Mc^{\U^*(\l)}$.
\end{thm}
\begin{proof}
By Corollary \ref{crlVerma_DCF} and Proposition \ref{thDCF_tw_pr},
there exists a dynamical twist in the category ${}_{\B\tl}\!\bar \O$,
which is a collection of morphisms ${}_{X}\!(F_{V,W}) \in \End_\l(X\tp V\tp
W)$,
where $V$ and $W$ are $\g$-modules and $X$ is a generic $\l$-module.
Using the natural filtration in generalized Verma modules,
one can prove that morphisms $\{{}_{X}\!(F_{V,W})\}$ are invertible in
$\bar \Mc^{\U^*(\l)}$.
The morphisms ${}_{X}\!(F_{V,W})$ defines a collection of $\l$-equivariant
maps $X^*\tp X\to \End_k(V\tp W)$, which
gives rise to a collection of  maps $ X^*\tp X\to \U(\g)\tp \U(\g)$
for generic $X$, since the dynamical twist is natural with respect to the
arguments $V$ and $W$.
This collection determines an $\l$-equivariant  map $ \bar \F\colon
\U^*(\l)\to \U(\g)\tp \U(\g)$
defined for generic elements $\U^*(\l)$. By Remark \ref{iso},
the dynamical category ${}_{\B\tl}\!\bar \O$ is isomorphic to  $\bar
\Mc^{\U^*(\l)}$.
Under this isomorphism, the dynamical twist $\{{}_{X}\!(F_{V,W})\}$ goes
over
to the map $\la\mapsto \bar \F(\la)$ for $\la\in X^*\tp X$ and $X$ generic,
which reduces to
the twisting cocycle (\ref{twistH*}) in representations.
\end{proof}
\section{Dynamical associative algebras and quantum vector bundles}
\label{sQVB}
\subsection{Classical vector bundles}
Let $H$ be a Lie group and $P$ be a principal  $H$-bundle.
Denote by $\A=\A(P)$ the algebra of functions on $P$.
Let $V$ be a finite dimensional left $H$-module.
An associated vector bundle $V(M)$ on $M=P/H$ with the
fiber $V$ is defined as the coset space $(P\times V)/H$ by
the action $(p,v)\mapsto(ph,h^{-1}v)$, $(p,v) \in (P\times V)$,
$h\in H$. The global sections of
$V(M)$ are identified with the space
$\bigl(\A(P)\ot V\bigr)^H\simeq\Hom_H\bigl(V^*,\A(P)\bigr)$. Let us denote
by $\A^V$ the space
of global sections of $V(M)$.
When $V=k$, the trivial module, the space $\A^k$ is
canonically identified with the  subalgebra in $\A$ of $H$-invariant
functions; in other words, $\A^k=\A(M)$.
The tensor product of vector bundles corresponds to
the tensor product of sections, which is induced
by multiplication in $\A$: given $s_V\in \A^V$ and $s_W\in \A^W$
the section $s_W\tp s_V \in \A^{V\tp W}$ is
$$
(s_W\tp s_V)(w\tp v):=s_W(w)s_V(v), \quad w\tp v\in W^*\tp V^*\simeq (V\tp
W)^*.
$$
In particular, the tensor product of sections makes the space $\A^V$ a
two-sided module over $\A^k$.
\subsection{Quantum vector bundles}
\label{QVB}
Fix a Hopf algebra $\Ha$ over the ground ring $k$ and consider a dynamical
associative
 algebra  $\A$ in the category $\bar \Mc^{\Ha^*}$,
cf. Example \ref{algH*}.
We are going to introduce associated vector bundles over the
"non-commutative coset space"
corresponding to the action of $\Ha$ on $\A$.
\begin{definition}
Let $V$ be a right $\Ha$-module. \select{The associated vector bundle}
$\A^V$
with fiber $V$ is a space of all $\Ha$-equivariant maps (sections)
$s_V\colon V^*\to \A$.
\end{definition}
\noindent

Observe that the restriction of the dynamical multiplication in $\A$
to a group-like element $\la$  of  $\Ha^*$ defines
a bilinear operation $\st{\la}\colon\A\tp \A\to \A$, which is
$\Ha$-equivariant,
since group-like elements are invariant in $\Ha^*$.
Now we can define a product of sections.
Let $V$ and $W$ be two $\Ha$-modules.
Take  $s_V\colon V^*\to \A$ and $s_W\colon W^*\to \A$ to be sections of
$\A^{V}$ and $\A^{W}$.
Fix a group-like element  $\la\in \Ha^*$.
The map $s_W\st{\la} s_V\colon(V\tp W)^*\simeq W^*\tp V^*\to \A$,
\be
\label{product_of_sections}
(s_W\st{\la} s_V)(w\tp v):= s_W(w)\st{\la} s_V(v), \quad w\tp v \in  W^*\tp
V^*,
\ee
is a section of the bundle $\A^{V\tp W}$.
The subspace of $\Ha$-invariants $\A^k\subset\A$ is obviously closed
under $\st{\la}$ for every group-like element $\la\in \Ha^*$.
\begin{thm}
\label{th_bundle}
For any group-like element $\la\in \Ha^*$ the multiplication $\st{\la}$
provides
 $\A^k$ with the structure of an associative algebra, $\A^k_\la$,
and makes the space $\A^V$ a left $\A_\la^k$-module.
If $V=k_\al$, i.e. is the 1-dimensional representation of $\Ha$ defined by
the character $\alpha$,
then the line bundle $\A^V$ is also a right $\A^k_{\la\al^{-1}}$-module
with respect to  $\st{\la}$. For any $a\in \A^k_\la$, $s_V\in \A^V$, and
$s_W\in  \A^W$
\be
a\st{\la}(s_V\st{\la}s_W)=(a\st{\la}s_V)\st{\la}s_W.
\ee
\end{thm}
\begin{proof}
Sections of the line bundle $\A^{k_\al}$ may be treated as elements  $a\in
\A$ such that $h\tr a = \al^{-1}(h) a$
for all $h\in \Ha$ (the inverse is understood in the sense of the algebra
$\Ha^*$).
For $a\in A^{k_\al}$ and $b,c\in\A$, the formula (\ref{multiplications})
turns into
\be
\label{group-like}
(a\st{\la}b) \st{\la}c=a\st{\la}(b\st{\la\al^{-1}}c)
\ee
under the assumption that $\la\in \Ha^*$ is group-like.
Setting $\al=1$ (the unit of $\Ha^*$) in (\ref{group-like}),  we find that
$\st{\la}$ is associative
when restricted to $\A^k$,
and makes  it an associative algebra, $\A_\la^k$. Also we see that
$\A$ is a left $\A_\la^k$-module.
This induces the structure of a left  $\A_\la^k$-module on $\A^V$ for
every
$\Ha$-module $V$, by formula (\ref{product_of_sections}).
Assuming $b,c\in \A^k$ in (\ref{group-like}) we obtain a right
$\A^k_{\la\al^{-1}}$-module
structure on the space $\A^{k_\al}$.
\end{proof}
\section{Vector bundles on semisimple coadjoint orbits}
\label{sVBCO}
The problem of equivariant quantization of function algebras on
semisimple coadjoint orbits of simple Lie groups was studied, e.g.,
in \cite{DGS1,DolJ,DM1,DM2,DM4,DS}. Quantization of vector bundles
on semisimple orbits as modules over the quantized function
algebras was considered in \cite{D1,GLS}.
In this section we apply dynamical associative algebras
to quantize the entire "algebra" of sections of all
vector bundles on semisimple coadjoint orbits.
\subsection{Dynamical quantization of the function algebra on a group}
Let $\g$ be a simple Lie algebra and $G$ the corresponding connected
Lie group.
We will apply the previous considerations to the problem
of equivariant quantization of vector bundles
on semisimple orbits in $\g^*$ with respect to  the coadjoint action of
$G$.
Denote by $\A(G)$ the algebra of polynomial functions on $G$.
The group $G$ acts on itself by the left and the right regular actions.
These actions induce two left commuting actions of $\U(\g)$
on $\A(G)$ via
the differential operators $\rho_1(x)$ and $\rho_2(x)$, $x\in\U(\g)$,
respectively . Here
$\rho_1(x)$ ($\rho_2(x)$) is
the differential operator on $G$ that is
the right (left) invariant extension of $\gamma(x)$ ($x$), where $\gamma$
denotes the antipode in $\U(\g)$.

Given a representation $\pi:G\to \End(V)$, one assigns to each
$f\in \End(V)^*$ the function $f\circ\pi$ on $G$. Identifying
$\End(V)^*$ with $\End(V)=V\ot V^*$ via the trace pairing,
these assignments give
the well known natural isomorphism
\be
\label{isonatfund}
\oplus_E E\ot E^*\to\A(G),
\ee
where $E$ runs
over all irreducible representations of $G$. Then the $\rho_1$-action
can be treated as an action on the $E$-factor while the $\rho_2$-action --
on the $E^*$-factor of each term $E\ot E^*$ in the direct sum
(\ref{isonatfund}).

As a manifold, a semisimple orbit is the quotient $M=G/H$, where $H$
is a Levi subgroup with the Lie algebra $\l\subset \g$.
Recall the decomposition $\l=\l_0\oplus\c$, where $\l_0$ is the semisimple part
and $\c$ the center of $\l$.
The manifold
$G$ may be considered as a principal $H$-bundle on $M$. Any
equivariant vector bundle on $M$ is associated to the
bundle $G$ via a representation $V$ of $H$. We denote this
vector bundle by $V(M)$.
It has the vector space $V$ as the fiber.

The global sections of the bundle $V(M)$ are identified with
the space $\bigl(\A(G)\ot V\bigr)^\l=\Hom_\l\bigl( V^*,\A(G)\bigr)$,
where $\A(G)$ is considered as a $\U(\g)$-module with respect
to the $\rho_2$-action. Since $M$ is an affine variety,
one can identify the vector bundle $V(M)$ with its global sections
and consider it as a $\U(\g)$-module with respect
to the $\rho_1$-action. In particular, suppose $V=\C_\lambda$ is
the one dimensional representation of $\l$ induced by the
character $\la\in \c^*$. Let us consider $\la$ as an element
of $\h^*$ via the natural embedding $\c^*\subset \h^*$.
Due to isomorphism (\ref{isonatfund}),
the assignment
$\Hom_\l\bigl(\C_\la^*,\A(G)\bigr)\ni\ff\mapsto \ff(1)\in\A(G)$
gives a natural isomorphism of $\U(\g)\ot\U(\l)$-modules
\be
\label{isonatf}
\Hom_\l\bigl(\C_\la^*,\A(G)\bigr)\to\A(G)[-\la],
\ee
where $\A(G)[-\la]$ is a subspace of $\l_0$-invariant elements of $\A(G)$
of weight $-\la$
with respect to the $\rho_2$-action. It is obvious that $\A(G)[-\la]$
is a $\U(\g)$-module with respect to the $\rho_1$-action and it is
naturally
isomorphic to
$\oplus_E E\ot E^*[-\la]$,
where $E^*[-\la]$ denotes the subspace of $\l_0$-invariant elements
of $E^*$ of weight $-\la$.
It is clear that isomorphism (\ref{isonatf}) is actually non-zero
only if $\la$ is an integer weight. In this case the map (\ref{isonatf})
identifies
$\A(G)[-\la]$ with the space of global sections of the line bundle
$\C_\la(M)$.
In particular, the function algebra on $M$ is naturally isomorphic to
the $\U(\g)$-module algebra
$\oplus_E E\ot E^*[0]\subset\A(G)$.

Applying the dynamical twist $\bar \F$ constructed in Theorem \ref{H*Utw}
to the $\U(\g)$-module algebra $\A(G)$ with respect to the
$\rho_2$-action,
we obtain a dynamical associative $\U(\l)$-algebra in the category $\Bar\O^{\U^*(\l)}$.
This algebra is equal to $\A(G)$ as a $\U(\g)$-module (with respect to
$\rho_1$-action)
and has
the family of multiplications parameterized by generic $\la\in\U^*(\l)$
and defined as $\bar m_\la=m\circ\bar \F^\la$,
where $m$ is the original multiplication in the algebra $\A(G)$.
Applying Theorem \ref{th_bundle} to this dynamical associative algebra, we
obtain a quantization
of vector bundles on $G/H$. Obviously, this quantization is equivariant
with respect to the $\rho_1$-action of $\U(\g)$.

Let us consider the dynamical twist $\bar \F^\la$ restricted to $\U^*(\c)$.
Applied to the subalgebra $\A(G)^{\l_0}\subset \A(G)$ of $\l_0$-invariant functions
on $G$ with respect to the $\rho_2$-action,
this restriction makes $\A(G)^{\l_0}$ a dynamical associative $\U(\l)$-algebra
over the base $\U^*(\c)$. As a $\U(\g)$-module, it is formed by sections of all linear bundles
on $M$. Let us describe this dynamical algebra in more detail.

Let $M_\la=\U(\g)\ot_{\U(\p^+)}\C_\la$ be the Verma module
corresponding to the one dimensional representation of $\U(\l)$
associated to $\la\in\c^*$. Let $\Hom^0(M_\la,M_\mu)$ denote
the subspace of locally finite elements in $\Hom_\C(M_\la,M_\mu)$
with respect to the adjoint action of $\U(\g)$. In fact,
$\Hom^0(M_\la,M_\mu)$ is not zero only when $\mu-\la$ is an
integer weight. Remind that for a $\U(\g)$-module $E$, we denote by
$E[\la]$, $\la\in\c^*$, the subspace of $\l_0$-invariant elements of weight~$\la$.

\begin{propn}\label{propnlf}
Let $V$ be a finite dimensional representation of $\U(\g)$.
Then, there is a natural morphism of $\U(\g)$-modules:
$V\ot V^*[\la-\mu]\to \Hom^0(M_\la,M_\mu)$,
$\la,\mu\in\c^*$, $\mu$ is generic. When $V$ is irreducible,
this morphism is embedding.
These embeddings give rise to the natural isomorphism
\be
\label{isonat}
j_{\la,\mu}:\oplus_E E\ot E^*[\la-\mu]\to \Hom^0(M_\la,M_\mu)
\ee
of $\U(\g)$-modules, where $E$ runs over all finite dimensional
irreducible representations of $\U(\g)$.
\end{propn}

\begin{proof} It is enough to prove the first part of the proposition
and to show that the multiplicity of $V$ in $\Hom^0(M_\la,M_\mu)$ is equal
to $\dim V^*[\la-\mu]$. Applying the Frobenius reciprocity,
one proves that for generic $\mu\in\c^*$ the space
$\Hom_{\U(\g)}(M_\la,M_\mu\ot V^*)$ is naturally
isomorphic to $V^*[\la-\mu]$; the proof is the same as in \cite{ES1}.
But $\Hom_{\U(\g)}(M_\la,M_\mu\ot V^*)\cong
\Hom_{\U(\g)}\bigl(V,\Hom(M_\la,M_\mu)\bigr)$,
which proves the proposition.
\end{proof}

Compositions $M_\nu\to M_\mu\to M_\la$ generate the map
$\Hom^0(M_\mu,M_\la)\ot \Hom^0(M_\nu,M_\mu)\to
\Hom^0\bigl(M_\nu,M_\la\bigr)$,
$\la,\mu,\nu\in\c^*$. Due to isomorphisms (\ref{isonat}) and
(\ref{isonatfund}), this map defines the morphism of
$\U(\g)\ot\U(\l)$-modules
\be
\label{multint}
\A(G)[\mu-\la]\ot\A(G)[\nu-\mu]\to\A(G)[\nu-\la].
\ee
Since $\A(G)[\beta]=0$ unless $\beta$ is an integer weight,
this morphism is defined for generic $\la\in\c^*$, i.e. for
$\la\not\in\mathcal Y$,
where $\mathcal Y$ is from Proposition \ref{Shap}.
Indeed,
if $\la\not\in\mathcal Y$, then also $\mu,\nu\not\in\mathcal Y$ when
the differences $\mu-\la$ and $\nu-\mu$ are integer weights.

Fixing a generic $\la$ in (\ref{multint}) and varying $\mu$ and $\nu$,
we obtain from (\ref{multint}) the morphism
\be
\label{multint11}
\A(G)^{\l_0}\ot\A(G)^{\l_0}\to\A(G)^{\l_0}.
\ee
These morphisms form a family of multiplications parameterized by elements
$e^\la\in\U^*(\c)$ (or $\la\in\c^*$) for generic $\la$.
Since the elements $e^\la$ form a basis of $\U^*(\c)$ over $\C$, this
family extends by linearity to all
generic elements of $ \U^*(\c)$.
One can check that this family makes $\A(G)^{\l_0}$ a dynamical associative
$\U(\l)$-algebra over the base $\U^*(\c)$.
Comparing the construction of this multiplication and the construction of
twist from Theorem \ref{H*Utw},
we come to the following.
\begin{propn}
\label{propninter}
The dynamical associative  multiplication (\ref{multint11})
has the form  $\bar m_\la=m\circ\bar\F^\la$ for generic $\la\in\c^*$,
where $\bar\F^\la$ is a dynamical twist over the base $\U^*(\l)$ from Theorem
\ref{H*Utw}.
\end{propn}

\subsection{Deformation quantization of the Kirillov brackets and
vector bundles on coadjoint orbits}
Dynamical twist from Theorem \ref{H*Utw} applied to $\A(G)$
gives a dynamical algebra, which, by Theorem \ref{th_bundle},
defines quantization of vector bundles on $M=G/H$
as left modules over the quantized algebra of functions on $M$.
This quantization of vector bundles is obviously $\g$-equivariant. Restricted
to the function algebra on $M$, it gives quantization of the Kirillov
brackets in the following way.

Let $t$ be an independent variable. Denote by $\g_t$ the Lie algebra
over $\C[[t]]$ with bracket $[x,y]_t:=t[x,y]$ for $x,y\in\g$, where
$[.,.]$ is the original bracket in $\g$. Then there is an algebra morphism
$\ff_t:\U(\g_t)\to \U(\g)[[t]]$ induced by the correspondence
$x\mapsto tx$ for $x\in\g$.
As was shown in \cite{DGS2}, the equivariant quantization
of the Kirillov Poisson bracket corresponding to
the semisimple orbit passing through $\la\in\c^*\subset\g^*$ is
identified with the image of $\U(\g_t)$ by the composition map
$\U(\g_t)\to\U(\g)[[t]]\to \Hom^0(M_{\la/t},M_{\la/t})$, where the first
map
is $\ff_t$ and the second one is the action map.

Using this fact, one can show
that the multiplication $\bar m_{\la,t}:=\bar m_{\la/t}$ from Proposition \ref{propninter}
being restricted to
$\A(M)=\A(G)[0]$ gives a
$\U(\g)$-equivariant deformation
quantization
of the Kirillov Poisson bracket when $M$ is realized as a coadjoint
orbit passing through $\la$.
Since the multiplication $\bar m_{\la,t}$
depends, in fact, on $\la/t$, we do not need that $\la$ to be generic in
the
deformation quantization.
Note that also for any formal path $\la(t)=\la_0+t\la_1+...\in\c^*[[t]]$
the multiplication
$\bar m_{\la(t),t}$ gives a $\U(\g)$-equivariant deformation
quantization on the orbit passing through $\la_0$, with the
appropriate Kirillov bracket.
\begin{remark} Any equivariant deformation quantization of
the Kirillov bracket on the orbit can be obtained in this way, and
different paths in $\c^*$ give not equivalent quantizations, \cite{D1}.
\end{remark}
For $\la\in\c^*$, consider $\bar m_\la=m\circ\bar\F^\la$, where $\bar\F^\la$ is the dynamical twist
from Theorem \ref{H*Utw} and $m$ the usual multiplication, as
a map $\A(G)^{\ot 2}\to \A(G)$.
In the similar way, one can prove
\begin{propn}\label{propndefq}
For any $\la_0\in\c^*$ and for any formal path
$\la(t)=\la_0+t\la_1+...\in\c^*[[t]]$,
the multiplication $\bar m_{\la(t),t}:=\bar m_{\la(t)/t}$ gives a
$\U(\g)\ot\U(\l)$-equivariant map $\A(G)^{\ot 2}\to \A(G)[[t]]$
that coincides modulo $t$ with the original
multiplication in $\A(G)$.
\end{propn}

Now, we fix $\la_0\in\c^*$ and consider $M=G/H$ as the semisimple
orbit passing through $\la_0$. Remind that the space of global sections
of the vector bundle
$V(M)$ corresponding to a representation $V$ of $H$ is identified with
the space $\bigl(\A(G)\ot V\bigr)^\l=\Hom_\l\bigl(V^*,\A(G)\bigr)$. We
identify $V(M)$
with its global sections.
Applying Theorem \ref{th_bundle} to the dynamical algebra
obtained from $\A(G)$ by the dynamical twist from Theorem \ref{H*Utw}
and using Proposition \ref{propndefq}, we obtain
\begin{thm}
\label{thmvecbund}
Let $\la(t)=\la_0+t\la_1+...$
be a formal path in $\c^*$. Then the dynamical multiplication
$\bar m_{\la(t)/t}$ defines a $\U(\g)$-equivariant multiplication
$\st{\la}$
on the global sections of equivariant vector bundles on $M$. This
multiplication
is a deformation of the usual tensor product of the sections
and satisfies the following properties:

1) Restricted to $\A(M)$, the operation $\st{\la(t)}$ defines a
deformation
quantization $\A_{\la(t)}(M)$
of the function algebra $\A(M)$ corresponding to the Kirillov bracket;

2) Let $s_V$ and $s_W$ are global section of vector bundles $V(M)$
and $W(M)$, and $a$ is a function on $M$. Then
$$a\st{\la}(s_V\st{\la}s_W)=(a\st{\la}s_V)\st{\la}s_W.$$
In particular, any vector bundle $V(M)$ is a left $\A_{\la(t)}(M)$-module
with respect to the action map $a\ot s\mapsto a\st{\la(t)}s$,
where $a\in\A_{\la(t)}(M)$, $s\in V(M)$.

3) The line bundle $\C_\alpha(M)$, where $\alpha\in\c^*$
is an integer weight, is also a right module over the algebra
$\A_{\la(t)-t\alpha}(M)$ with respect
to the action map
$s\ot b\mapsto s\st{\la(t)}b$, where
$s\in\C_\alpha(M)$, $b\in\A_{\la(t)-t\alpha}(M)$.
\end{thm}

\subsection{The quantum group case}
Let $\U_q(\g)$ be the Drinfeld-Jimbo quantum group corresponding to $\g$
and
$\U_q(\l)$ be considered as its quantum subgroup corresponding to
the Levi subalgebra $\l\subset\g$. Let $\A_q(G)$ denote the dual
algebra to $\U_q(\g)$ consisting of matrix elements of finite
dimensional representations of $\U_q(\g)$.
The algebra $\A_q(G)$ is a quantization of the classical
algebra $\A(G)$, it is equivariant under the left and the right regular
actions of $\U_q(\g)$,
which we replace by two left actions, $\rho_1$ and $\rho_2$, as above.

Let $\bar \F_q$ be the dynamical twist constructed in Theorem
\ref{H*Utw}
with help of generalized Verma modules over $\U_q(\g)$.
Applying this twist to the algebra $\A_q(G)$,
we obtain a dynamical algebra $\bigl(\A(G),\bar m_{q,\la}\bigr)$
in the category $\Bar\O^{\U_q^*(\l)}$.
This algebra is equal to $\A_q(G)$ as a $\U_q(\g)$-module
(with respect to $\rho_1$-action)
and has
the family of multiplications $\bar m_{q,\la}$ parameterized by
generic $\la\in\U_q^*(\l)$
and defined by  $\bar m_{q,\la}:=m_q\circ \bar \F_{q,\la}$,
where $m_q$ is the original multiplication in  $\A_q(G)$.
It is obvious that $\A_{q,\la}(G)$ is a $\U_q(\g)$-module algebra
with respect to the $\rho_1$-action.

One can show that replacing simultaneously $\la$ by $\la/t$ and $q$ by
$q^t$,
we obtain the family of multiplications $\bar m_{q^t,\la,t}:=\bar
m_{q^t,\la/t}$
that gives a $\U_{q^t}(\g)$-equivariant deformation quantization of
$\A(M)$.
Also, there is a $q$-analog of Theorem \ref{thmvecbund} which reduces to
Theorem \ref{thmvecbund} when $q=1$.


\begin{thebibliography}{A}
\bibitem[ABB]{ABB} J. Avan, O. Babelon, E. Billey:
 {\em The Gervais-Neveu-Felder equation and the quantum Calogero-Moser
systems},
 Commun. Math. Phys.
 {\bf 178} (1996) 281 -- 299.
\bibitem[AL]{AL} A. Alekseev, A. Lachowska:
 {\em Invariant $*$-product on coadjoint orbits and the Shapovalov
pairing},
arXiv: math.QA/0308100.
\bibitem[AM]{AM} A. Alekseev, E. Meinrenken:
 {\em The non-commutative Weil algebra},
 Invent. Math.
 {\bf 135} (2000) 135--172.
\bibitem[AF]{AF} A. Alekseev, L. Faddeev:
{\em $T^*(G)_t$: a toy model of conformal field theory},
Commun. Math. Phys.
 {\bf 141} (1991) 413--422.
\bibitem[BDFh]{BDFh} J. Balog, L. Dabrowski, and L. Feh\'er:
{\em Classical $r$-matrix and exchange algebra in WZNW and Toda field
theories},
Phys. Lett. B,
{\bf 244} \# 2 (1990) 227--234.
\bibitem[Dr1]{Dr1} V.  Drinfeld:
  {\em  Quantum Groups}, in Proc. Int. Congress of Mathematicians,
  Berkeley,  1986,  ed. A. V. Gleason, AMS, Providence (1987) 798--820.
\bibitem[Dr2]{Dr2} V.  Drinfeld:
{\em Hamiltonian structures on Lie groups, Lie bialgebras and the geometric
meaning
of the classical Yang-Baxter euqations},
Sov. Math. Dokl.
{\bf 27} (1983) 68--71.
\bibitem[Dr3]{Dr3} V.  Drinfeld:
  {\em  Quasi-Hopf algebras},
  Leningrad Math.J.
  {\bf 1} (1990) 1419--1457.
\bibitem[Dr4]{Dr4} V.  Drinfeld:
  {\em  On Poisson homogeneous spaces of Poisson-Lie groups},
  Theor. Math. Phys., {\bf 95} (1993) 226-227.
\bibitem[DolJ]{DolJ} B. P. Dolan, O. Jahn,
{\em Fuzzy Complex Grassmannian Spaces and their Star Products},
hep-th/0111020.
\bibitem[D1]{D1} J. Donin:
{\em $\U_h(\g)$-invariant quantization of coadjoint orbits and
vector bundles over them},
 J. Geom. Phys.,
 {\bf 38}  \#1 (2001) 54--80.
\bibitem[D2]{D2} J. Donin:
 {\em Quantum $G$-manifolds},
 Contemp. Math., {\bf 315}  (2002) 47--60.
\bibitem[DGS1]{DGS1} J. Donin, D. Gurevich, S. Shnider:
 {\em Double quantization in some
  orbits in the coadjoint representiations of simple Lie groups},
  Commun. Math. Phys.
  {\bf 204} (1999) 39--60.
\bibitem[DGS2]{DGS2} J. Donin, D. Gurevich, S. Shnider:
 {\em Quantization of function algebras on semisimple orbits in $\g^*$},
arXiv:q-alg/9607008.
\bibitem[DM1]{DM1} J. Donin, A. Mudrov:
{\em Method of quantum characters in equivariant quantization},
Commun. Math. Phys.
{\bf 234} (2003)  533--555.
\bibitem[DM2]{DM2} J. Donin, A. Mudrov:
{\em Explicit  equivariant quantization on coadjoint orbits of  GL(n)},
Lett. Math. Phys.,
{\bf 62} (2002) 17--32.
\bibitem[DM3]{DM3} J. Donin and A. Mudrov: {\em Reflection Equation, Twist,
and Equivariant Quantization},
 Isr. J. Math.,
 {\bf 136} (2003) 11--28.
\bibitem[DM4]{DM4} J. Donin, A. Mudrov:
{\em  Quantum coadjoint orbits of GL(n) and generalized Verma modules};
math.QA/0212318.
\bibitem[DM5]{DM5} J. Donin and A. Mudrov: {\em Quantum groupoids
associated with dynamical categories},
arXiv: math.QA/0311116.
\bibitem[DO]{DO} J. Donin and V. Ostapenko:
{\em Equivariant quantization on quotients of simple Lie groups by
reductive subgroups of maximal rank},
Czech. Journ. of Physics,
{\bf 52} \# 11 (2002), 1213--1218.
\bibitem[DS]{DS} J. Donin and S. Shnider:
{\em Quantum symmetric spaces},
J. Pure \& Applied Alg.
{\bf 100} (1995) 103--116.
\bibitem[EE1]{EE1} B. Enriquez, P. Etingof:
 {\em Quantization of Alekseev-Meinrenken dynamical r-matrices},
arXiv:math.QA/0302067.
\bibitem[EE2]{EE2} B. Enriquez, P. Etingof:
{Quantization of classical dynamical r-matrices withnonabelian base},
arXiv:math.QA/0311224.
\bibitem[EK]{EK} P. Etingof,  D. Kazhdan:
    {\em Quantization of Lie bialgebras},
    Selecta Math. {\bf 2}, \# 1  (1996) 1--41.
\bibitem[ES1]{ES1} P. Etingof, O. Schiffmann:
 {\em Lectures on the dynamical Yang-Baxter equation. Quantum Groups and
Lie theory},
London Math. Soc. Lecture Note {\bf 290} (2001) 89--129.
\bibitem[ES2]{ES2} P. Etingof, O. Schiffmann:
 {\em On the moduli space of classical dynamical r-matrices},
 Math. Res. Lett.
 {\bf 9} (2001)  157--170.
\bibitem[ESS]{ESS} P. Etingof, O. Schiffmann, T. Schedler:
 {\em Explicit quantization of dynamical r-matrices for finite dimensional
semisimple Lie algebras},
J.AMS
 {\bf 13} (2000)  595--609.
\bibitem[EV1]{EV1} P. Etingof, A. Varchenko:
 {\em Geometry and classification of solutions of the classical dynamical
Yang-Baxter equation},
Commun. Math. Phys.
 {\bf 192} (1998)  77--120.
\bibitem[EV2]{EV2} P. Etingof, A. Varchenko:
 {\em Solutions of the quantum dynamical Yang-Baxter equation and dynamical
quantum groups},
Commun. Math. Phys.
 {\bf 196} (1998)  591--640.
\bibitem[EV3]{EV3} P. Etingof, A. Varchenko:
 {\em Exchange dynamical quantum groups},
Commun. Math. Phys.
 {\bf 205} (1999)  19--52.
\bibitem[F]{F} G. Felder:
 {\em Conformal field theories and integrable models associated to elliptic
curves},
 Proc. ICM Zurich,
 {\bf } (1994)  1247--1255.
\bibitem[Fh]{Fh} L. Feh\'er:
{\em Dynamical $r$-matrices and Poisson-Lie symmetries of the chiral WZNW
model},
arXiv:hep-th/02012006.
\bibitem[FhMrsh]{FhMrsh} L. Feh\'er and I. Marshall:
{\em On a Poisson-Lie analogue of the classical dynamical Yang-Baxter
equation for self dual Lie algebras},
arXiv:math.QA/0208159.
\bibitem[Fad]{Fad} L. Faddeev:
{\em On the exchande matrix of the WZNW model},
 Commun. Math. Phys.
{\bf 132} (1990) 131--138.
\bibitem[FRT]{FRT} L. Faddeev,  N. Reshetikhin, and L. Takhtajan:
 { \em Quantization of Lie groups and Lie algebras},
 Leningrad Math. J. {\bf 1}  (1990)  193--226.
\bibitem[GLS]{GLS} D. Gurevich, R. Leclercq, P. Saponov:
 {\em q-Index on braided spheres}, math.QA/0207268.
\bibitem[GN]{GN} J.-L. Gervais, A. Neveu:
 {\em Novel triangle relation and absence of tachyons in Liouville string
theory},
 Nucl. Phys. B
 {\bf 238} (1984)  125--141.
\bibitem[J]{J} J.-C. Jantzen:
{\em Kontravariante Formen auf  induzierten Darstellungen halbeinfacher aa
Lie-Algebren},
Math. Ann. {\bf 226} (1977) 53--65.
\bibitem[Kas]{Kas} Ch. Kassel:
{\em Quantum groups}, Springer, NY 1995.
\bibitem[Kar]{Kar} E. Karolinskii:
  {\em A classification of Poisson homogeneous spaces of complex reductive
Poisson-Lie groups},
   In: Banach Center Publ.
  {\bf 51} (1996) 103--108.
\bibitem[KMST]{KMST} E. Karolinsky, K. Muzykin, A. Stolin, V. Tarasov:
{\em Dynamical Yang-Baxter equations, quasi-Poisson homogeneous spaces, and
quantization},
arXiv:math.QA/0309203.
\bibitem[KSkl]{KSkl} P. P. Kulish, E. K. Sklyanin:
  {\em Algebraic  structure  related to the reflection equation,}
  J.  Phys.  A  {\bf  25}  (1992)  5963--5976.
\bibitem[KS]{KS}  P.  P.  Kulish,  R.  Sasaki:
  {\em Covariance properties of reflection  equation  algebras,}
  Prog.  Theor. Phys. {\bf 89} $\# 3$  (1993)  741--761.
\bibitem[Lu1]{Lu1} J.-H. Lu:
  {\em  Classical dynamical r-matrices and homogeneous Poisson structures
on $G/H$ and $K/T$},
  Commun. Math. Phys.
  {\bf 212} (2000)  337--370.
\bibitem[Lu2]{Lu2} J. H. Lu:
 {\em Hopf algebroids and quantum gouppoids},
  Int. J. Math., {\bf 7} (1996) 47--70.
\bibitem[O]{O} V. Ostrik,
  {\em Module categories, weak Hopf algebras, and modular invariants},
  arXiv: math.QA/111139.
\bibitem[RS]{RS}  N. Reshetikhin, M. Semenov-Tian-Shansky:
{\em Quantum $R$-matrices and factorization problem},
  J. Geom. Phys. {\bf 5} (1988), 533--550.
\bibitem[S]{S} O. Schiffmann:
 {\em On classification of dynamical r-matrices},
Math. Res. Lett.
 {\bf 5} (1998) 13--30.
\bibitem[Sem]{Sem}  M. Semenov-Tian-Shansky: {\em Poisson-Lie Groups,
Quantum Duality Principle, and the Quantum Double},
  Contemp. Math. {\bf 175} (1994) 219--248.
\bibitem[Xu1]{Xu1} P. Xu:
 {\em Triangular dynamical r-matrices and quantization},
Adv. Math.
 {\bf 166} (2002)  1--49.
\bibitem[Xu2]{Xu2} P. Xu:
 {\em Quantum Dynamical Yang-Baxter Equation Over a Nonabelian  basis},
Commun. Math. Phys.
 {\bf 226} (2002)  475--495.
\end{thebibliography}
\end{document}